\documentclass[journal]{IEEEtran}
\usepackage{newunicodechar}
\usepackage[table]{xcolor}
\newunicodechar{₂}{$_2$}
\usepackage{amsmath}
\allowdisplaybreaks
\makeatletter
\def\section{\@startsection{section}{1}{\z@}{1.2ex plus 0.3ex minus 0.3ex}{0.4ex plus 0.2ex}{\normalfont\normalsize\centering\scshape}}
\def\subsection{\@startsection{subsection}{2}{\z@}{1.0ex plus 0.2ex minus 0.2ex}{0.2ex}{\normalfont\normalsize\itshape}}
\makeatother
\newcommand{\keepdisplaytogether}{\interdisplaylinepenalty=10000 \predisplaypenalty=0 }
\usepackage[official]{eurosym}
\usepackage{amssymb}
\usepackage{stfloats}
\usepackage{amsfonts}
\usepackage{nomencl}
\usepackage{algorithm} \usepackage{algpseudocode}
\usepackage[hidelinks]{hyperref}
\usepackage{alphalph}
\makeatletter
\renewcommand*{\@alph}[1]{\alphalph{#1}}
\makeatother
\usepackage{amsthm}
\usepackage{upgreek}
\newtheorem{proposition}{Proposition}

\theoremstyle{remark}

\usepackage{etoolbox}
\usepackage{booktabs}
\usepackage{cite}
\usepackage{mathtools}
\usepackage{cancel}
\usepackage{multirow}
\usepackage{graphicx}
\graphicspath{{../figures/}{figures/}{./}}
\usepackage{tikz}
\newcommand{\Phase}[1]{\Statex \textbf{#1}}

\usetikzlibrary{arrows.meta,positioning,calc,backgrounds,fit}
\makenomenclature
\renewcommand\nomgroup[1]{%
  \item[\bfseries
  \ifstrequal{#1}{P}{Parameters}{%
  \ifstrequal{#1}{V}{Variables}{%
  \ifstrequal{#1}{I}{Indices and Sets}{}}}%
]}
\AtBeginDocument{%
  \setlength{\abovedisplayskip}{2pt plus 1pt minus 1pt}%
  \setlength{\belowdisplayskip}{2pt plus 1pt minus 1pt}%
  \setlength{\abovedisplayshortskip}{1pt plus 1pt}%
  \setlength{\belowdisplayshortskip}{1pt plus 1pt}%
}
\selectfont
\begin{document}
\setlength{\textfloatsep}{4pt plus 2pt minus 2pt}
\setlength{\dbltextfloatsep}{4pt plus 2pt minus 2pt}
\setlength{\floatsep}{4pt plus 1pt minus 1pt}
\setlength{\dblfloatsep}{4pt plus 1pt minus 1pt}
\setlength{\intextsep}{4pt plus 1pt minus 1pt}
\setlength{\abovecaptionskip}{1pt plus 1pt}
\setlength{\belowcaptionskip}{0pt}
\title{Incentive-Based Regulatory Mechanism for Tariff Design in District Heating and Cooling Networks}
\author{%
Hossein Hesamzadeh, Wicak Ananduta, Eric Dumont, Zacharie De Greve%
}
\maketitle
\begin{abstract}
  Low-temperature district heating and cooling networks exchange heat bidirectionally at near-ambient temperatures. When one supplier owns the shared loop, it controls both the tariff and the network temperature that sets every user’s conversion efficiency. Existing price-cap and cost-based rules leave this non-price instrument unregulated and reach neither the first-best dispatch nor an explicit division of gains. This paper proposes a two-part, incentive-based tariff that anchors the supplier’s lump-sum fee, per building and hour, to an audited outside alternative and divides the participation gains through one regulated share. We model a mixed-integer bilevel Stackelberg game in which the supplier sets hourly tariffs and the temperature regime, anticipating the buildings’ co-optimization of electricity and thermal exchange. We prove that the mechanism attains the first-best welfare, posts the shadow price of the efficient dispatch as its tariff, keeps every building at or above its outside option every hour, and keeps the supplier solvent whenever operation adds surplus. A Belgian case study confirms these properties. The mechanism reproduces the centralized welfare whenever operation is worthwhile, declines to operate otherwise, holds the loop at its first-best temperature, and outperforms cap-based regulation and Ramsey pricing in welfare, user protection, and cost recovery.
\end{abstract}
\section*{Nomenclature}

{\small\renewcommand{\IEEEiedlistdecl}{\setlength{\itemsep}{0pt}\setlength{\parsep}{0pt}\setlength{\topsep}{1pt}}%
\noindent\textbf{\textit{Indices and Sets}}
\begin{IEEEdescription}[\IEEEsetlabelwidth{$\Delta\mathrm T^{\mathrm{diff}},\Delta\mathrm T^{\mathrm{int}},\Delta\mathrm T^{\min}$}]
\item[$t,\ \mathcal T$] Hour index; horizon of length $T$.
\item[$b,\ \mathcal B$] Building index; set of $B$ buildings.
\item[$k,\ \mathcal K$] Temperature-bin index and set.
\item[$g,\ \mathcal G$] Device index; set $\{\mathrm{GSHP},\mathrm{GSCh}\}$.
\item[$\mathcal K^{\mathrm{cold}}_{b}$] Bins allowing direct cooling for $b$.
\end{IEEEdescription}
\noindent\textbf{\textit{Parameters}}
\begin{IEEEdescription}[\IEEEsetlabelwidth{$\Delta\mathrm T^{\mathrm{diff}},\Delta\mathrm T^{\mathrm{int}},\Delta\mathrm T^{\min}$}]
\item[$\uplambda^{\mathrm{El,imp}}_t,\uplambda^{\mathrm{El,exp}}_t,\uplambda^{\mathrm{gas}}$] Grid import, export; gas prices, \euro/kWh.
\item[$\mathrm{COP}^{\mathrm{Hp}}_{k},\mathrm{COP}^{\mathrm{Ch}}_{k}$] Building HP and chiller COP per bin.
\item[$\mathrm{COP}^{g}_{k},\overline{\mathrm P}^{g}$] Device $g$ COP per bin; capacity, kW.
\item[$\overline{\mathrm P}^{\mathrm{Hp}}_{b},\overline{\mathrm P}^{\mathrm{Ch}}_{b}$] Building HP, chiller caps, kW.
\item[$\hat{\mathrm H}^{h}_{b,t},\hat{\mathrm H}^{c}_{b,t}$] Reference heat, cold demand, kWh$_\mathrm{th}$.
\item[$\widehat{\mathrm H}^{\mathrm{DRC}}_{b},\hat{\mathrm P}^{\mathrm{PV}}_{b,t}$] Direct-cooling rating; PV forecast.
\item[$\mathrm d^{\mathrm h}_{b},\mathrm d^{\mathrm c}_{b}$] Comfort curvatures, \euro/kWh$_\mathrm{th}^{2}$.
\item[$\upvarepsilon,\bar{\mathrm p},\bar{\mathrm H}$] Elasticity; mean unit cost, hourly load.
\item[$\upgamma^{\mathrm{pump}},\upgamma^{\mathrm{DRC}}$] Well-pump (kW per m$^3$/h), DRC (kW$_\mathrm{e}$/kW$_\mathrm{th}$) power ratios.
\item[$\upalpha_{w}$] Wellhead transfer coeff., kWh$_\mathrm{th}$/m$^3$.
\item[$\uprho^{w},\mathrm c^{w},\mathrm{kA}^{\mathrm{netw}}$] Water density, heat capacity; pipe loss, kW/K.
\item[$\mathrm V^{\mathrm{netw}},\mathrm V^{\mathrm{ACC}},\overline{\mathrm S}^{\mathrm{ACC}}$] Network, accumulator volumes; capacity.
\item[$\Delta\mathrm T^{\mathrm{diff}},\Delta\mathrm T^{\mathrm{int}},\Delta\mathrm T^{\min}$] Pipe separation, bin interval, pinch, K.
\item[$\mathrm T^{\mathrm{low}}_{k},\mathrm T^{\mathrm{up}}_{k},\mathrm T^{\min}_{\mathrm{netw,c}}$] Bin edges; cold-pipe minimum, \textdegree C.
\item[$\mathrm T^{\mathrm{soil}}_{t},\mathrm T^{\mathrm{ww}},\mathrm T^{\mathrm{cw}}$] Soil, warm-, cold-well temperatures, \textdegree C.
\item[$\mathrm T^{\mathrm{C,sup}}_{b},\mathrm T^{\mathrm{C,ret}}_{b}$] Cooling supply, return temperatures, \textdegree C.
\item[$\mathrm s,\mathrm s_{\max}$] Surplus share; adequate maximum.
\item[$\mathrm{h}_d,\Delta\mathrm{t},\mathrm M_{T}$] Day heating share; time step; big-$M$ (K).
\item[$\upeta,\mathrm{EF}$] Fee-norm weight; CO$_2$ factor, kg/kWh.
\item[$\mathrm R,\overline{\uplambda}^{H},\overline{\uplambda}^{C}$] S1 headroom; heat, cold cap references.
\item[$\mathrm{COP}^{\mathrm{ASHP}}_{t},\mathrm{COP}^{\mathrm{ASCh}}_{t}$] Outside-option ASHP, ASCh COP.
\item[$\upeta^{\mathrm{GB}}$] Outside-option boiler efficiency.
\item[$\overline{\mathrm P}^{\mathrm{ASHP}}_{b},\overline{\mathrm P}^{\mathrm{GB}}_{b},\overline{\mathrm P}^{\mathrm{ASCh}}_{b}$] Outside-option capacities, kW.
\end{IEEEdescription}
\noindent\textbf{\textit{Upper-level variables, $\Omega^{\mathrm{UL}}$}}
\begin{IEEEdescription}[\IEEEsetlabelwidth{$\Delta\mathrm T^{\mathrm{diff}},\Delta\mathrm T^{\mathrm{int}},\Delta\mathrm T^{\min}$}]
\item[$\lambda^{\mathrm{Th}}_{t},\theta_{k,t}$] Thermal tariff (\euro/kWh); bin selector.
\item[$T^{\mathrm{netw,w}}_{t},T^{\mathrm{netw,c}}_{t}$] Warm-, cold-pipe temperatures, \textdegree C.
\item[$T^{\mathrm{netw,M}}_{t},\Delta T_{t}$] Mean temperature (\textdegree C); warming step (K).
\item[$T^{\theta}_{k,t},P^{\theta,g}_{k,t}$] Bin mean temperature; supplier bin power.
\item[$P^{\mathrm{Sup}}_{t}$] Supplier electrical draw, kWh$_\mathrm{e}$.
\item[$P^{\mathrm{Sup},g}_{t},H^{\mathrm{Sup},g}_{t}$] Device $g$ draw (kWh$_\mathrm{e}$), output (kWh$_\mathrm{th}$).
\item[$H^{\mathrm{loss,wp}}_{t},H^{\mathrm{loss,cp}}_{t}$] Warm-, cold-pipe losses, kWh$_\mathrm{th}$.
\item[$H^{\mathrm{netw}}_{t}$] Network inertia, kWh$_\mathrm{th}$.
\item[$u^{\mathrm H}_{t},u^{\mathrm C}_{t}$] Warm-, cold-well flows, m$^{3}$.
\item[$S^{\mathrm{ACC}}_{t},H^{\mathrm{ACC}}_{t}$] Accumulator state, kWh$_\mathrm{th}$; charge, kW$_\mathrm{th}$.
\item[$y^{\mathrm{inc}}_{t},y^{\mathrm{dec}}_{t}$] Network warming, cooling binaries.
\end{IEEEdescription}
\noindent\textbf{\textit{Lower-level variables, $\Omega^{\mathrm{LL}}$}}
\begin{IEEEdescription}[\IEEEsetlabelwidth{$\Delta\mathrm T^{\mathrm{diff}},\Delta\mathrm T^{\mathrm{int}},\Delta\mathrm T^{\min}$}]
\item[$P^{\theta,\mathrm{Hp}}_{b,k,t},P^{\theta,\mathrm{Ch}}_{b,k,t}$] Building HP, chiller bin powers, kWh$_\mathrm{e}$.
\item[$H^{\mathrm{Th,B}}_{b,t}$] Metered net thermal exchange, kWh$_\mathrm{th}$.
\item[$H^{\mathrm{DRC}}_{b,t}$] Direct-cooling flow, kWh$_\mathrm{th}$.
\item[$P^{\mathrm{El,imp}}_{b,t},P^{\mathrm{El,exp}}_{b,t}$] Grid import; PV export, kWh$_\mathrm{e}$.
\end{IEEEdescription}
\noindent\textbf{\textit{Other quantities}}
\begin{IEEEdescription}[\IEEEsetlabelwidth{$\Delta\mathrm T^{\mathrm{diff}},\Delta\mathrm T^{\mathrm{int}},\Delta\mathrm T^{\min}$}]
\item[$x^{\mathrm h}_{b,t},x^{\mathrm c}_{b,t}$] Served heat, cool levels, kWh$_\mathrm{th}$.
\item[$U_{b,t}$] Consumption utility, \euro.
\item[$\Phi_{b,t},\rho_b$] Regulatory fee; rebate, \euro.
\item[$D_{b,t},\hat D_{b,t}$] Hourly participation gain; announced, \euro.
\item[$G_b$] Horizon participation gain of $b$, \euro.
\item[$\Pi^{\mathrm{Sup}}_{t}$] Supplier hourly profit, \euro.
\item[$\Pi^{\mathrm B}_{b,t},\Pi^{\mathrm{Alt}}_{b,t}$] Building, outside-option payoffs, \euro.
\item[$\mathrm C^{\mathrm B}_{b,t},\mathrm C^{\mathrm{Alt}}_{b,t}$] On-, off-network purchase costs, \euro.
\item[$\lambda^{*}_{t}$] Shadow price of exchanged heat, \euro/kWh.
\item[$W_t,W,\overline W,\Delta W$] Hourly, total, maximum welfare; incremental surplus, \euro.
\item[$P^{\mathrm{ASHP}}_{b,t},P^{\mathrm{GB}}_{b,t},P^{\mathrm{ASCh}}_{b,t}$] Outside-option device inputs, kWh.
\item[$P^{\mathrm{a,El}}_{b,t},P^{\mathrm{a,s}}_{b,t},P^{\mathrm{a,exp}}_{b,t}$] Outside-option grid draw; PV use, export.
\end{IEEEdescription}
}
\section{Introduction}\label{sec:intro}
Decarbonizing European Union (EU) energy sectors places District Heating (DH) at the center of low-carbon urban supply, especially when coupled with renewable electricity, thermal sources, and storage. DH operates under regulated or deregulated schemes \cite{LI201556}, but its high upfront network cost makes it a natural monopoly, so its supplier holds discretion over both tariffs and the operating conditions that shape users' feasible responses \cite{doi/10.2833/962525}. The instruments that constrain this discretion range from ex-post competition-law oversight to dedicated sectoral regulation \cite{BILLERBECK2023113377}. Low-Temperature District Heating and Cooling (LTDHC) sharpens the problem, since lower temperatures, tighter electricity coupling, bidirectional exchange, and distributed renewables give the supplier more ways to shape what users can do \cite{BUFFA2019504}. Competition in EU heating markets is weak \cite{de2025study}, so market pressure disciplines neither the tariff nor the operating conditions. This discipline must come from sector-specific pricing rules.

In the literature, local pricing is approached from two perspectives. Centralized cooperative schemes let the connected agents jointly optimize and then divide the joint gain, through Nash-bargaining energy sharing \cite{Cui21}, dynamic net metering in energy communities \cite{alahmed}, or peer-to-peer multi-energy trading \cite{Jin25}. Non-cooperative schemes instead let a profit-maximizing leader set prices and operating conditions while the users respond, either without regulatory bounds \cite{Fan18,Che19b} or under them \cite{Lu21,Jin21,LI20,Bai22}. Cooperation presumes an agreed objective, shared cost information, and enforceable side payments. None of these exists for a profit-maximizing LTDHC supplier with private cost data, so the non-cooperative view is the one that matches the institutional setting and is adopted here. These studies differ in what constrains the supplier's price. They are reviewed in turn, together with work on conventional DH networks, where local thermal pricing is more developed.

\textit{Non-Cooperative Deregulated Pricing (market-based).}
The first stream unbundles the network from supply and lets a market clear on it, so that the heat price emerges as the dual of the network-constrained dispatch, as with locational marginal heat prices under combined-heat-and-power Stackelberg clearing \cite{Fan18} and an optimal thermal-flow market beneath a strategic provider \cite{Che19b}. Such prices measure the marginal value of heat correctly, but presuppose organized clearing, third-party network access, and a separation of network ownership from supply, which DH rarely has. Entry barriers and owner control persist \cite{BILLERBECK2023113377} and deregulated systems are confined to a few EU markets \cite{PAZERAITE2022112888}, so efficient prices must be induced by regulation, the question this stream leaves open.

\textit{Non-Cooperative Regulated Pricing.}
The second stream sets supplier tariffs under exogenous bounds, such as a cap tied to the upstream energy price \cite{Lu21}, \cite{Jin21}, a unified electricity--thermal retailer above responsive buildings \cite{LI20}, and a capped low-temperature tariff mobilizing building flexibility \cite{Bai22}. A cap limits how much the supplier may charge, but not how it distributes charges across hours and exchange directions, so allocative efficiency is never a design objective, and the welfare forgone remains unquantified. Two further limitations bind in LTDHC. First, any cap needs a defensible external reference price, yet the marginal value of bidirectional near-ambient exchange changes hour by hour and depends on the loop temperature the supplier chooses, so a reference fixed in advance cannot track it. Second, capping the tariff does not exhaust the supplier's withholding power, because raising or lowering the loop temperature changes every device's efficiency and decides which users may cool directly from the loop, and those non-price instruments stay unregulated. Another classical mechanism, regulator-set rates under a break-even constraint (Ramsey pricing, e.g., \cite{alahmed}), maximizes welfare subject to the firm covering its cost, which leaves the utility no margin on days when that constraint binds. Neither instrument reaches the first-best (the dispatch that maximizes total welfare) and neither constrains the loop temperature.

A thermal supplier and an electricity transmission operator share a similar structure, with network-bound users and a natural monopoly whose efficiency must be induced by regulation, not competition \cite{SODERHOLM2011742}. There, the binding constraint is information asymmetry, since the regulator observes neither the firm's costs nor the users' valuations, and the classic menu-and-audit response of optimal-regulation theory \cite{baron1982regulating,laffont1993theory} is informationally demanding. Loeb and Magat \cite{loeb1979decentralized} showed that a subsidy equal to full consumer surplus yields first-best pricing without cost knowledge, but it is publicly funded and merely shifts the information asymmetry to the demand side. The Incremental Surplus Subsidy \cite{sappington1988regulating} shrinks the Loeb--Magat transfer to the period-over-period welfare gain, yet keeps subsidy financing and a surplus-measurement burden. On the price-cap side, the Vogelsang--Finsinger process \cite{vogelsang1979regulatory} steers a multiproduct monopoly toward Ramsey prices through lagged tariff constraints, Vogelsang \cite{vogelsang1989twopart} carries the same adjustment into the fixed part of a two-part tariff, so that the users rather than a subsidy finance the incentive, and \cite{hogan2010merchant} applies that tariff to transmission. The Hesamzadeh--Rosell\'on--Gabriel--Vogelsang (HRGV) mechanism \cite{hesamzadeh2018simple} merges the surplus-subsidy and price-cap traditions, evolving the fee by the realized change in total surplus and converging to the first-best over successive regulatory periods. Its transplant from transmission to thermal networks is not immediate. The HRGV mechanism is dynamic, anchored to a base period the firm can distort. It imposes no individual-rationality constraint on the paying users, so a user's fee may exceed its benefit and drive it off the network. Implementations remain investment-oriented bilevel programs \cite{khastieva2020transmission}, and in the variant of \cite{xia2024incentive} surplus sharing is obtained only through a fee-tuning parameter that distorts efficiency.
To the knowledge of the authors, no existing incentive mechanism for a local thermal monopoly, which has no public subsidy base and no organized market to lean on, is subsidy-free, settled within a single horizon without a manipulable base period, participation-guaranteeing, and equipped with a surplus-sharing parameter at no welfare cost.

This paper closes that gap. The proposed mechanism, built on the HRGV approach, keeps the incremental-surplus principle but anchors the fee to an audited outside alternative rather than to the firm's own history, which removes the base period the firm can distort. Because that alternative serves the same reference demand under the same contracted comfort terms, the settlement is carried by submetered flows, posted prices, and the certified reference technology, and no valuation or demand curve needs to be obtained from users by the regulator. The transfer is financed by the served users rather than by a subsidy. The bilevel structure extends the regulator's reach to the network temperature, which price-only instruments leave free. A single regulated share divides the resulting gains without moving the dispatch. Taken together, our main contributions are as follows:
\begin{enumerate}
    \item \textit{Incentive-based and welfare-aligned local pricing policy.} The mechanism recasts the HRGV incremental-surplus fee logic \cite{hesamzadeh2018simple} into a two-part tariff for a price-making monopoly thermal utility. Unlike HRGV, where the fee follows the aggregate surplus variation across regulatory periods, it defines a building-level fee relative to an auditable outside alternative and a regulated share of the participation gains.
    \item \textit{Temperature-aware bilevel monopoly model.} Making the supplier's seasonal temperature control explicit, we reformulate the operational model of \cite{wirtz2021temperature} as a bilevel Stackelberg game in which the monopolist jointly sets thermal prices and the network temperature while anticipating the buildings' co-optimization of electricity and bidirectional exchange, unlike price-only bilevel pricing \cite{Lu21,Jin21} and locational marginal heat pricing \cite{Fan18,Che19b}.
    \item \textit{Proofs of welfare alignment under monopoly.} We prove that the mechanism is efficient and leaves every building at or above its outside option in every hour, together with its best-response, settlement, cost-recovery, and pricing properties, under an explicit timing protocol and information requirements. We verify the properties numerically on an operating Belgian LTDHC network.
\end{enumerate}
The remainder of the paper is organized as follows. Section~\ref{sec:proposed_model} presents the bilevel model and the proposed mechanism. Section~\ref{sec:case_study} benchmarks the mechanism against a centralized welfare optimum, a price-cap regulation, and Ramsey pricing. Section~\ref{sec:conclusion} concludes. The proofs of the theoretical developments are in Appendix~\ref{app:proofs} and the outside-alternative benchmark in Appendix~\ref{app:alt}.

\section{Mathematical Modeling} \label{sec:proposed_model}
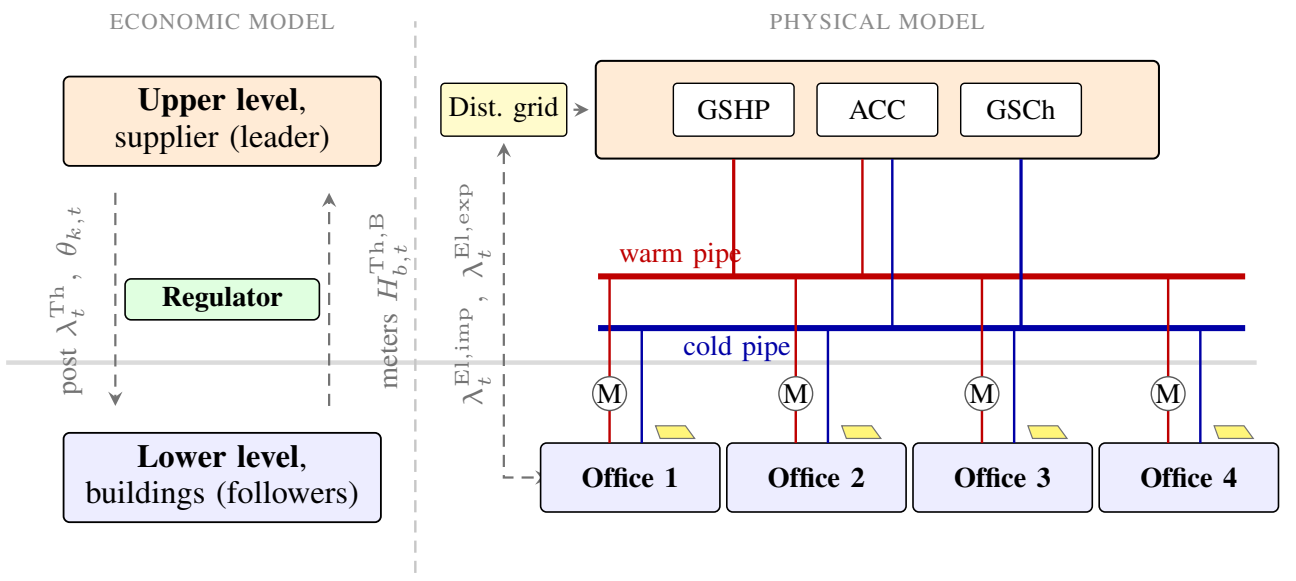
\begin{figure}[t]
\centering
\resizebox{\linewidth}{!}{%
\begin{tikzpicture}[
  font=\footnotesize, >=latex,
  bldg/.style={draw, rounded corners=1.5pt, fill=blue!7, minimum width=1.32cm, minimum height=0.60cm, align=center, inner sep=1pt, line width=0.4pt},
  plnt/.style={draw, rounded corners=1.5pt, fill=orange!16, minimum width=4.90cm, minimum height=0.85cm, align=center, inner sep=2.5pt, line width=0.5pt},
  sub/.style={draw, rounded corners=1pt, fill=white, minimum height=0.46cm, minimum width=1.05cm, align=center, inner sep=2pt, font=\scriptsize, line width=0.35pt},
  econ/.style={draw, rounded corners=1.5pt, text width=2.55cm, align=center, inner sep=3pt, line width=0.45pt},
  warm/.style={red!75!black, line width=1.5pt},
  cold/.style={blue!65!black, line width=1.5pt},
  sig/.style={->, >=stealth, line width=0.55pt, densely dashed, black!55, shorten >=1.5pt, shorten <=1.5pt},
  stp/.style={circle, draw=black!60, fill=white, inner sep=0.4pt, minimum size=8.5pt, font=\tiny, line width=0.4pt},
  met/.style={circle, draw=black!65, fill=white, inner sep=0.2pt, minimum size=8pt, font=\scriptsize, line width=0.35pt},
]
\path (-0.45,-1.05) (11.20,3.95);
\draw[black!14, line width=1.1pt] (0.02,0.85) -- (10.90,0.85);
\draw[densely dashed, black!22, line width=0.6pt] (3.58,-0.98) -- (3.58,3.82);
\node[font=\scriptsize, text=black!45] at (7.60,3.82) {\textsc{physical model}};
\node[font=\scriptsize, text=black!45] at (1.90,3.82) {\textsc{economic model}};
\node[plnt] (P) at (7.60,3.05) {};
\node[sub] (gshp) at (6.35,3.05) {GSHP};
\node[sub] (acc)  at (7.60,3.05) {ACC};
\node[sub] (gsch) at (8.85,3.05) {GSCh};
\node[sub, fill=yellow!25] (ret) at (4.35,3.05) {Dist.\ grid};
\draw[sig] (ret.east) -- (P.west);
\draw[sig, <->] (ret.south) -- (4.35,-0.15) -- (4.80,-0.15);
\node[font=\scriptsize, text=black!55, rotate=90, anchor=center] at (4.08,1.45) {$\uplambda^{\mathrm{El,imp}}_t,\ \uplambda^{\mathrm{El,exp}}_t$};
\draw[warm] (5.17,1.60) -- (10.80,1.60);
\draw[cold] (5.17,1.15) -- (10.80,1.15);
\node[warm, anchor=west, font=\scriptsize] at (5.21,1.78) {warm pipe};
\node[cold, anchor=west, font=\scriptsize] at (5.77,0.97) {cold pipe};
\draw[warm, line width=0.8pt] (6.35,2.62) -- (6.35,1.60);
\draw[cold, line width=0.8pt] (8.85,2.62) -- (8.85,1.15);
\draw[warm, line width=0.6pt] (7.47,2.62) -- (7.47,1.60);
\draw[cold, line width=0.6pt] (7.73,2.62) -- (7.73,1.15);
\foreach \i/\nm in {1/{Office 1},2/{Office 2},3/{Office 3},4/{Office 4}}{
  \pgfmathsetmacro\bx{5.45 + (\i-1)*1.62}
  \node[bldg, minimum width=1.56cm] (b\i) at (\bx,-0.15) {\scriptsize\textbf{\nm}};
  \draw[fill=yellow!65, draw=black!55, line width=0.3pt]
    (\bx+0.26,0.19) -- (\bx+0.56,0.19) -- (\bx+0.48,0.30) -- (\bx+0.22,0.30) -- cycle;
  \draw[warm, line width=0.6pt] (\bx-0.18,0.15) -- (\bx-0.18,1.60);
  \draw[cold, line width=0.6pt] (\bx+0.10,0.15) -- (\bx+0.10,1.15);
  \node[met] at (\bx-0.18,0.58) {M};
}
\node[econ, fill=orange!16] (UL) at (1.90,2.95) {\textbf{Upper level}, supplier (leader)};
\node[econ, fill=blue!7]    (LL) at (1.90,-0.15) {\textbf{Lower level}, buildings (followers)};
\draw[sig] (0.97,2.38) -- (0.97,0.42);
\node[font=\scriptsize, text=black!55, rotate=90, anchor=center] at (0.55,1.40) {post $\lambda^{\mathrm{Th}}_t,\ \theta_{k,t}$};
\draw[sig] (2.83,0.42) -- (2.83,2.38);
\node[font=\scriptsize, text=black!55, rotate=90, anchor=center] at (3.36,1.40) {meters $H^{\mathrm{Th,B}}_{b,t}$};
\node[draw, rounded corners=1.5pt, fill=green!12, text width=1.55cm, align=center, inner sep=2pt, font=\scriptsize, line width=0.4pt] (R) at (1.90,1.40) {\textbf{Regulator}};
\end{tikzpicture}}
\caption{The LTDHC system and its three players. The supplier posts tariff and temperature, buildings self-dispatch, and the regulator sets the fee rule~\eqref{eq:reg-body} and the rebate~\eqref{eq:rho-proportional} and audits the settlement. The physical model illustrates the
case studied in Section \ref{sec:case_study}.}
\label{fig:network}
\end{figure}

As Fig.~\ref{fig:network} shows, the setting involves three parties sharing a single, local, bidirectional LTDHC loop, a profit-maximizing monopolistic supplier that owns the network and sets its prices, the connected buildings that act as thermal prosumers, and a regulatory body that administers a pricing mechanism. The economic model in the left panel of Fig.~\ref{fig:network} captures this hierarchy. The supplier commits first, posting the hourly tariff and the temperature policy, and each building can only react to those posted terms. The supplier is therefore the leader of the resulting Stackelberg game and the buildings are its followers, with the regulator auditing the settlement between them. The physical model in the right panel, which illustrates how the assets of all players are connected, is the loop this game is played on. Every building carries settlement-grade meters that record its electricity use, its photovoltaic (PV) generation, and its two-way heat exchange with the loop. Hour by hour, a building's heat pumps and chillers either draw heat from the loop or return surplus heat or cold to it, depending on the temperature the supplier holds in the loop. Each actor's model is defined in turn below.
\subsection{Upper-level (Thermal supplier problem)}
\label{sec:upper_level}
This subsection defines the supplier problem. The supplier owns the plant, the pipes, and the accumulator tank, sets the network temperature as the leader, and sees the buildings only through their metered exchange. Its objective function is
\begin{align}
& \max_{\Omega^{\mathrm{UL}}}\ \sum_{t\in\mathcal T}\!\Big(\Pi_t^{\mathrm{Sup}}+\!\sum_{b\in\mathcal B}\!\Phi_{b,t}\Big), \label{eq:UL-obj-short}\\
& \Pi_t^{\mathrm{Sup}}\!=\!\lambda_t^{\mathrm{Th}}\!\!\sum_{b\in\mathcal B}\!\!H_{b,t}^{\mathrm{Th,B}}\!-\!\uplambda_t^{\mathrm{El,imp}}P_t^{\mathrm{Sup}},
& \forall t,\label{eq:UL-profit-ates-pen}
\end{align}
where $\Pi^{\mathrm{Sup}}_t$ is the supplier's hourly thermal-sales profit net of electricity cost, with $\lambda^{\mathrm{Th}}_t$ being the announced thermal price and $P^{\mathrm{Sup}}_t$ being the supplier's grid purchase at retail price $\uplambda^{\mathrm{El,imp}}_t$. The second term of \eqref{eq:UL-obj-short} is the fee revenue, where $\Phi_{b,t}$ is the regulatory fee building $b$ pays in hour $t$, defined by the rule \eqref{eq:reg-body} of Section~\ref{sec:mechanism}, and $\Omega^{\mathrm{UL}}$ collects the upper-level variables. The constraints below, modeling the thermal network, follow the LTDHC formulation of~\cite{wirtz2021temperature}:
\setlength{\jot}{1pt}
\begin{subequations}\label{UL}
\begin{flalign}
& P_{t}^{\mathrm{Sup}} =
  P^{\mathrm{Sup,GSHP}}_{t} + P^{\mathrm{Sup,GSCh}}_{t}
  + \upgamma^{\mathrm{pump}}(u^{\mathrm{H}}_{t}+u^{\mathrm{C}}_{t}),
  && \forall t, \label{eq:UL-a}
\end{flalign}
\begin{flalign}
& H_{t}^{\mathrm{Sup,GSHP}} - H_{t}^{\mathrm{Sup,GSCh}} = \nonumber\\
& \quad \sum_{b \in \mathcal{B}} H^{\mathrm{Th,B}}_{b,t}
  + H_{t}^{\mathrm{loss,wp}}
  - H_{t}^{\mathrm{loss,cp}}
  + H_{t}^{\mathrm{netw}}
  + H_{t}^{\mathrm{ACC}},
  && \forall t, \label{eq:UL-b}
\end{flalign}
\begin{flalign}
& H_{t}^{\mathrm{netw}} =
  \frac{\uprho^{w}\,\bigl(\mathrm{V}^{\mathrm{netw}}+\mathrm{V}^{\mathrm{ACC}}\bigr)\,\mathrm{c}^{w}}{3600\,\Delta\mathrm{t}}
  \,\Delta T_{t},
  && \forall t, \label{eq:UL-c}
\end{flalign}
\begin{flalign}
& H_{t}^{\mathrm{loss,wp}} = \mathrm{kA}^{\mathrm{netw}}
  \bigl(T_{t}^{\mathrm{netw,w}} - \mathrm{T}_{t}^{\mathrm{soil}}\bigr),
&& \forall t, \label{eq:UL-d}\\
& H_{t}^{\mathrm{loss,cp}} = \mathrm{kA}^{\mathrm{netw}}
  \bigl(\mathrm{T}_{t}^{\mathrm{soil}} - T_{t}^{\mathrm{netw,c}}\bigr),
&& \forall t, \label{eq:UL-d2}
\end{flalign}
\begin{flalign}
& \textstyle\sum_{k}\theta_{k,t} = 1,\quad T_{t}^{\mathrm{netw,M}} = \sum_{k} T^{\theta}_{k,t},
  && \forall t, \label{eq:UL-e}
\end{flalign}
\begin{flalign}
& \theta_{k,t}\in\{0,1\},\quad \theta_{k,t}\,\mathrm{T}^{\mathrm{low}}_{k}
  \le T^{\theta}_{k,t}
  \le \theta_{k,t}\,\mathrm{T}^{\mathrm{up}}_{k},
  && \forall k,t, \label{eq:UL-f}
\end{flalign}
\begin{flalign}
& \mathrm{T}^{\mathrm{up}}_{k} =
  \mathrm{T}^{\min}_{\mathrm{netw,c}} + \tfrac{\Delta \mathrm{T}^{\mathrm{diff}}}{2} + k\,\Delta \mathrm{T}^{\mathrm{int}},
&& \forall k, \label{eq:UL-h}\\
& \mathrm{T}^{\mathrm{low}}_{k} = \mathrm{T}^{\mathrm{up}}_{k} - \Delta \mathrm{T}^{\mathrm{int}},
&& \forall k, \label{eq:UL-h2}
\end{flalign}
\begin{flalign}
& T_{t}^{\mathrm{netw,w}} = T_{t}^{\mathrm{netw,M}} + \tfrac{\Delta \mathrm{T}^{\mathrm{diff}}}{2} = T_{t}^{\mathrm{netw,c}} + \Delta \mathrm{T}^{\mathrm{diff}},
&& \forall t, \label{eq:UL-i}
\end{flalign}
\begin{flalign}
& T_{t}^{\mathrm{netw,w}} - T_{t-1}^{\mathrm{netw,w}} = \Delta T_{t},
  && \forall t, \label{eq:UL-j}\\
& S_{t}^{\mathrm{ACC}} = S_{t-1}^{\mathrm{ACC}} + \Delta\mathrm{t}\, H_{t}^{\mathrm{ACC}},
  && \forall t, \label{eq:UL-l}
\end{flalign}
\begin{flalign}
& S_{0}^{\mathrm{ACC}} := S_{T}^{\mathrm{ACC}},\quad T_{0}^{\mathrm{netw,w}} := T_{T}^{\mathrm{netw,w}}, &&
  \label{eq:UL-m}
\end{flalign}
\begin{flalign}
& 2\,\overline{\mathrm{S}}^{\mathrm{ACC}}y_{t}^{\mathrm{inc}}
  \le S_{t}^{\mathrm{ACC}} + S_{t-1}^{\mathrm{ACC}}
  \le 2\,\overline{\mathrm{S}}^{\mathrm{ACC}}(1 - y_{t}^{\mathrm{dec}}),
  && \forall t, \label{eq:UL-ACC-yinc}\\
& y_{t}^{\mathrm{inc}}+y_{t}^{\mathrm{dec}}\le 1,\quad y_{t}^{\mathrm{inc}},y_{t}^{\mathrm{dec}}\in\{0,1\},
  && \forall t, \label{eq:UL-ACC-y}
\end{flalign}
\begin{flalign}
& -y_{t}^{\mathrm{dec}}\,\mathrm{M}_{T} \le \Delta T_{t} \le y_{t}^{\mathrm{inc}}\,\mathrm{M}_{T},
  && \forall t, \label{eq:UL-n}
\end{flalign}
\begin{flalign}
& H^{\mathrm{Sup},g}_{t} = \textstyle\sum_{k}\mathrm{COP}^{g}_{k}\,P^{\theta,g}_{k,t},\quad P^{\mathrm{Sup},g}_{t} = \sum_{k} P^{\theta,g}_{k,t},
&& \forall g,t, \label{eq:UL-p}
\end{flalign}
\begin{flalign}
& P^{\theta,g}_{k,t} \le \theta_{k,t}\,\overline{\mathrm{P}}^{g},
  && \forall g,k,t, \label{eq:UL-q}
\end{flalign}
\begin{flalign}
& H^{\mathrm{Sup,GSHP}}_{t} - P^{\mathrm{Sup,GSHP}}_{t} = \upalpha_{w}\,u^{\mathrm{H}}_{t},
&& \forall t, \label{eq:UL-s}\\
& H^{\mathrm{Sup,GSCh}}_{t} + P^{\mathrm{Sup,GSCh}}_{t} = \upalpha_{w}\,u^{\mathrm{C}}_{t},
  && \forall t, \label{eq:UL-t}
\end{flalign}
\end{subequations}
As given in \eqref{eq:UL-a}, the supplier's electrical draw is the sum of the compressor power of the Ground-Source Heat Pump (GSHP) and of the Ground-Source Chiller (GSCh) and the electricity of the well pumps. The thermal balance \eqref{eq:UL-b} allocates the supplier's net thermal output to the buildings' metered residuals, the pipe losses, the heat inertia of the network water, and the accumulator charging. The inertia term is defined in \eqref{eq:UL-c}, and \eqref{eq:UL-d} and \eqref{eq:UL-d2} give the pipe losses to the soil. The partition of unity in \eqref{eq:UL-e} keeps exactly one temperature bin active per hour and recovers the mean network temperature, \eqref{eq:UL-f} declares the binary selector and pins the mean temperature inside the active bin, whose edges are set by \eqref{eq:UL-h} and \eqref{eq:UL-h2}, and \eqref{eq:UL-i} fixes the warm and cold pipes at a constant separation around that mean. The warm pipe evolves by the temperature step \eqref{eq:UL-j} and the accumulator state by the recurrence \eqref{eq:UL-l}, both closing cyclically through \eqref{eq:UL-m}. Following \cite{wirtz2021temperature}, the network may warm (cool) only in hours when the accumulator is fully hot (cold). Constraints \eqref{eq:UL-ACC-yinc} and \eqref{eq:UL-ACC-y} enforce this with mutually exclusive direction binaries, and \eqref{eq:UL-n} ties the warming step to them. The thermal outputs of the two supplier devices, indexed by $g\in\mathcal G=\{\mathrm{GSHP},\mathrm{GSCh}\}$, follow bin-wise Coefficient of Performance (COP) curves in \eqref{eq:UL-p}, which also aggregates their electrical draws, each per-bin power is gated by the active bin in \eqref{eq:UL-q}, and the device source sides are linked to the well flows through \eqref{eq:UL-s} and \eqref{eq:UL-t}, where $\upalpha_{w}=\uprho^{w}\mathrm{c}^{w}(\mathrm{T}^{\mathrm{ww}}-\mathrm{T}^{\mathrm{cw}})/3600$ is the doublet wellhead energy-transfer coefficient. Direct Cooling (DRC) depends only on the loop temperature. Building $b$ can use it only when the active bin is cold enough, i.e.\ $k\in\mathcal{K}^{\mathrm{cold}}_{b}\subset\mathcal{K}$, the precomputed bins for which both pipes stay below the building's pinch temperatures, $T^{\mathrm{netw,c}}\!+\!\Delta\mathrm{T}^{\min}\!\le\!\mathrm{T}^{\mathrm{C,sup}}_{b}$ and $T^{\mathrm{netw,w}}\!+\!\Delta\mathrm{T}^{\min}\!\le\!\mathrm{T}^{\mathrm{C,ret}}_{b}$. Warming the loop out of these bins forces the building onto its compressor chiller.
\subsection{Lower-level (buildings problem)}\label{sec:building_constraints}
This subsection presents the optimization problem of the buildings. Each building operates three thermal assets behind its meter, a Water-Source Heat Pump (WSHP), a Water-Source Chiller (WSCh), and a direct-cooling connection, alongside a rooftop PV array, which is not a decision asset. Instead, its forecast output $\hat{\mathrm P}^{\mathrm{PV}}_{b,t}$ enters the electrical balance \eqref{eq:LL-Pbal} as a fixed parameter. Other dispatchable assets, such as batteries, can be added as long as they admit convex models. Each building sees only the posted price $\lambda^{\mathrm{Th}}_{t}$ and the network temperatures $T^{\mathrm{netw,w}}_{t},T^{\mathrm{netw,c}}_{t}$ implied by the upper-level bin policy $\theta_{k,t}$, and chooses how to serve its own demand. Its objective function is the net payoff
\begin{align}
& \max_{\Omega^{\mathrm{LL}}}\ \sum_{t\in\mathcal T}\bigl(\Pi^{\mathrm{B}}_{b,t}-\Phi_{b,t}\bigr), & \forall b,\label{eq:LL-obj}\\
& \Pi^{\mathrm{B}}_{b,t}\!=\!U_{b,t}\!-\!\bigl(\lambda_t^{\mathrm{Th}}H_{b,t}^{\mathrm{Th,B}}\!+\!\uplambda_t^{\mathrm{El,imp}}P_{b,t}^{\mathrm{El,imp}}\nonumber\\
& \qquad\qquad\!-\!\uplambda_t^{\mathrm{El,exp}}P_{b,t}^{\mathrm{El,exp}}\bigr), & \forall b,t,\label{eq:HRGV_Pi_B}\\
& U_{b,t}=
  -\mathrm{d}^{\mathrm{h}}_b\bigl(x^{\mathrm{h}}_{b,t}-\hat{\mathrm{H}}^{h}_{b,t}\bigr)^{2}\nonumber\\
& \qquad\qquad -\mathrm{d}^{\mathrm{c}}_b\bigl(x^{\mathrm{c}}_{b,t}-\hat{\mathrm{H}}^{c}_{b,t}\bigr)^{2}, & \forall b,t,\label{eq:HRGV_utility}\\
& x^{\mathrm{h}}_{b,t}\!=\!\textstyle\sum_{k}\mathrm{COP}^{\mathrm{Hp}}_{k}P^{\theta,\mathrm{Hp}}_{b,k,t},\nonumber\\
& x^{\mathrm{c}}_{b,t}\!=\!\textstyle\sum_{k}\mathrm{COP}^{\mathrm{Ch}}_{k}P^{\theta,\mathrm{Ch}}_{b,k,t}\!+\!H^{\mathrm{DRC}}_{b,t},
  & \forall b,t.\label{eq:HRGV_x}
\end{align}
The building payoff $\Pi^{\mathrm{B}}_{b,t}$ in \eqref{eq:HRGV_Pi_B} is the comfort utility $U_{b,t}$ minus the hourly costs of thermal exchange and grid electricity, net of PV exports at the grid export tariff $\uplambda^{\mathrm{El,exp}}_t$. Following a common demand-modeling practice \cite{jamal}, the utility \eqref{eq:HRGV_utility} is a quadratic deviation function, in which heating and cooling are penalized separately by the squared deviation of the supplied levels $x^{\mathrm{h}}_{b,t}$ and $x^{\mathrm{c}}_{b,t}$ in \eqref{eq:HRGV_x} from their respective targets. Here, the coefficients $\mathrm{d}^{\mathrm{h}}_{b}$ and $\mathrm{d}^{\mathrm{c}}_{b}$ govern price responsiveness. Direct cooling draws only the circulation-pump electricity $\upgamma^{\mathrm{DRC}}H^{\mathrm{DRC}}_{b,t}$. The constraints of the lower-level problem, with their dual variables in parentheses, are as follows:
\begin{subequations}\label{LL}
\begin{flalign}
& (\nu_{b,t})\!:\ H^{\mathrm{Th,B}}_{b,t}\!=\!\textstyle\sum_{k}\bigl(\mathrm{COP}^{\mathrm{Hp}}_{k}\!-\!1\bigr)P^{\theta,\mathrm{Hp}}_{b,k,t} \nonumber\\
& \qquad\;\; -\textstyle\sum_{k}\bigl(\mathrm{COP}^{\mathrm{Ch}}_{k}\!+\!1\bigr)P^{\theta,\mathrm{Ch}}_{b,k,t}\!-\!H^{\mathrm{DRC}}_{b,t},
  && \forall b,t, \label{eq:LL-res}
\end{flalign}
\begin{flalign}
& (\delta_{b,t})\!:\ P^{\mathrm{El,imp}}_{b,t}\!+\!\hat{\mathrm{P}}^{\mathrm{PV}}_{b,t}\!=\!\textstyle\sum_{k}\!\bigl(P^{\theta,\mathrm{Hp}}_{b,k,t}\!+\!P^{\theta,\mathrm{Ch}}_{b,k,t}\bigr) \nonumber\\
& \qquad\qquad\qquad\quad\;\;\, +\upgamma^{\mathrm{DRC}}H^{\mathrm{DRC}}_{b,t}\!+\!P^{\mathrm{El,exp}}_{b,t}, && \forall b,t, \label{eq:LL-Pbal}
\end{flalign}
\begin{flalign}
& (\phi^{\mathrm{El}}_{b,t})\!:\ P^{\mathrm{El,imp}}_{b,t}\!\ge\!0, && \forall b,t, \label{eq:LL-imp}\\
& (\underline{\phi}^{\mathrm{El,exp}}_{b,t},\overline{\phi}^{\mathrm{El,exp}}_{b,t})\!:\ 0\!\le\!P^{\mathrm{El,exp}}_{b,t}\!\le\!\hat{\mathrm{P}}^{\mathrm{PV}}_{b,t}, && \forall b,t, \label{eq:LL-exp}
\end{flalign}
\begin{flalign}
& (\nu^{\theta\mathrm{Hp}}_{b,k,t},\mu^{h}_{b,k,t})\!:\ 0\!\le\!P^{\theta,\mathrm{Hp}}_{b,k,t}\!\le\!\theta_{k,t}\overline{\mathrm{P}}^{\mathrm{Hp}}_{b}, && \forall b,k,t, \label{eq:LL-cap}\\
& (\nu^{\theta\mathrm{Ch}}_{b,k,t},\mu^{c}_{b,k,t})\!:\ 0\!\le\!P^{\theta,\mathrm{Ch}}_{b,k,t}\!\le\!\theta_{k,t}\overline{\mathrm{P}}^{\mathrm{Ch}}_{b},
  && \forall b,k,t, \label{eq:LL-cap2}
\end{flalign}
\begin{flalign}
& (\nu^{\mathrm{DRC}}_{b,t},\omega_{b,t})\!:\ 0\!\le\!H^{\mathrm{DRC}}_{b,t}\!\le\!\widehat{\mathrm{H}}^{\mathrm{DRC}}_{b}\!\sum_{k\in\mathcal{K}^{\mathrm{cold}}_{b}}\!\theta_{k,t},
  && \forall b,t. \label{eq:LL-DRC}
\end{flalign}
\end{subequations}
Constraint~\eqref{eq:LL-res} defines the building's thermal balance, where the metered net exchange with the loop equals the heat drawn by the heat pump, minus the heat rejected by the chiller and any direct cooling. Constraint~\eqref{eq:LL-Pbal} is the electrical balance. Grid import and on-site PV cover the compressors, the direct-cooling pump, and any PV export. Grid draw and export stay within the bounds \eqref{eq:LL-imp} and \eqref{eq:LL-exp}. Constraints~\eqref{eq:LL-cap} and \eqref{eq:LL-cap2} cap each device's power, which can be nonzero only in the active temperature bin, since exactly one bin is active in each hour. Constraint~\eqref{eq:LL-DRC} limits direct cooling to the heat-exchanger rating. It is available only while a cold bin is active, through the gate $\sum_{k\in\mathcal{K}^{\mathrm{cold}}_{b}}\theta_{k,t}$. The set $\Omega^{\mathrm{LL}}$ defines the lower-level variables, and the multipliers at the left of each constraint are the dual variables, $\Omega^{\mathrm{dual}}\!=\!\{\nu$, $\delta$, $\phi^{\mathrm{El}}$, $\underline{\phi}^{\mathrm{El,exp}}$, $\overline{\phi}^{\mathrm{El,exp}}$, $\nu^{\theta\mathrm{Hp}}_{k}$, $\mu^{h}_{k}$, $\nu^{\theta\mathrm{Ch}}_{k}$, $\mu^{c}_{k}$, $\nu^{\mathrm{DRC}}$, $\omega\}$.
\subsection{Tariff design}\label{sec:mechanism}
The proposed tariff has two parts. The volumetric part is the posted price $\lambda^{\mathrm{Th}}_{t}$, which steers the dispatch, and the lump-sum part is the regulatory fee, which transfers each building's participation gain to the supplier while leaving the building a rebate $\rho_b$. Financing the incentive transfer through the fixed part of a two-part tariff paid by the served users, rather than through public funds, follows Vogelsang \cite{vogelsang1989twopart}. His fixed part returns the firm's lagged profit and his price reaches marginal cost only as the adjustment converges, whereas here it is set from each building's gain over an audited outside alternative and the marginal price is posted at once. The rebate, in turn, treats the users' participation margin, which the HRGV mechanism leaves untreated. The regulatory body defines the regulatory fee $\Phi_{b,t}$ through the following feasibility problem:%
\begingroup\keepdisplaytogether
\begin{subequations}\label{eq:reg-body}
\begin{align}
&\rlap{\kern-2em$\text{Find }\Phi_{b,t},\ \forall b\!\in\!\mathcal B,\ t\!\in\!\mathcal T,\ \text{s.t.}$}\nonumber\\
& \Phi_{b,0}\!=\!D_{b,0}\!-\!\rho_b, & \forall b,\label{eq:HRGV_phi_init}\\
& \Phi_{b,t}\!=\!\Phi_{b,t-1}\!+\!D_{b,t}\!-\!D_{b,t-1}, & \forall b,\;t\!\ge\!1,\label{eq:HRGV_phi}\\
& D_{b,t}\!=\!\Pi^{\mathrm{B},*}_{b,t}(\lambda^{\mathrm{Th}}_{t},\theta_{k,t})\!-\!\Pi^{\mathrm{Alt}}_{b,t}, & \forall b,t,\label{eq:HRGV_D}
\end{align}
\end{subequations}\endgroup
where \eqref{eq:HRGV_phi_init} sets the fee at $t=0$, with $D_{b,0}=0$, \eqref{eq:HRGV_phi} updates it hour by hour by the change in the gain, so that the fee always equals the current gain net of the rebate, and \eqref{eq:HRGV_D} defines the gain $D_{b,t}$ as the payoff difference between being on and off the network. Here, $\Pi^{\mathrm{B},*}_{b,t}(\lambda^{\mathrm{Th}}_{t},\theta_{k,t})$ is the optimal payoff of building $b$, parameterized by the posted upper-level decisions, i.e., the objective \eqref{eq:HRGV_Pi_B} evaluated at the maximizer of the lower-level problem \eqref{eq:LL-obj}--\eqref{LL}. To lighten notation, the star and the arguments are dropped in what follows. The term $\Pi^{\mathrm{Alt}}_{b,t}$ is the hourly payoff of the building's audited outside alternative, the optimum of the small convex quadratic program of Appendix~\ref{app:alt} rather than a fixed reference bill, and solving it per building and hour from certified technology data makes the anchor auditable. $D_{b,t}$ is thus the building's hourly avoided-cost surplus from network membership.
\begin{proposition}[Best-Response Consistency]\label{prop:brc}
If the fee schedule $(\Phi^{*}_{b,t})$ is computed by \eqref{eq:reg-body} from the announced gains and held fixed at settlement, then, for any announced $(\lambda^{\mathrm{Th}}_t,\theta_{k,t})_{t\in\mathcal T}$, the prescribed lower-level response is the unique best response of every building, i.e., any unilateral deviation from it strictly lowers the building's net payoff.\end{proposition}
Proposition~\ref{prop:brc} states that the fee does not distort the buildings' response. Since the fee is a lump sum fixed at announcement, each building maximizes its payoff \eqref{eq:HRGV_Pi_B} at the posted price. The announced dispatch is therefore the one that takes place, and no building gains by deviating.
\begin{proposition}[Anonymous Settlement]\label{prop:info}
Let $\mathrm C^{\mathrm B}_{b,t}$ and $\mathrm C^{\mathrm{Alt}}_{b,t}$ be the purchase costs of building $b$ in hour $t$ on and off the network, and $\upvarepsilon$ the own-price elasticity to which $\mathrm d_b$ of \eqref{eq:HRGV_utility} is calibrated. Under the fee \eqref{eq:reg-body}, $\Phi_{b,t}$ is a function of metered purchases, posted prices, the certified reference technology, and the contracted comfort parameters only, and not of the supplier's cost function. Moreover, $D_{b,t}\to\mathrm C^{\mathrm{Alt}}_{b,t}-\mathrm C^{\mathrm B}_{b,t}$ as $\upvarepsilon\to0$, and $D_{b,t}-(\mathrm C^{\mathrm{Alt}}_{b,t}-\mathrm C^{\mathrm B}_{b,t})=\mathcal O(\upvarepsilon)$ if $\upvarepsilon>0$.
\end{proposition}
Proposition~\ref{prop:info} states that the fee can be settled from verifiable data alone. The regulator reads the meters, the posted and public prices, and the certified technology, and needs neither the supplier's cost function nor any user's valuation of comfort, which makes the mechanism anonymous in the sense of \cite{sibley1989asymmetric}. 
On the other hand, the sharing rule is a regulatory choice, applied to each building's gains aggregated over the settlement horizon, $\sum_t D_{b,t}$. Since \eqref{eq:HRGV_phi_init} and \eqref{eq:HRGV_phi} give $\Phi_{b,t}=D_{b,t}-\rho_b$ in every hour, summing the fee over the horizon shows that the building retains $\sum_t\Pi^{\mathrm{Alt}}_{b,t}+\mathrm{s}\sum_t D_{b,t}$ and the supplier $\sum_t\Pi^{\mathrm{Sup}}_t+(1-\mathrm{s})\sum_{b,t}D_{b,t}$ whenever the floor below is slack, with $\mathrm{s}$ the regulatory parameter:
\begin{equation}\label{eq:rho-proportional}
\rho_b=\max\Bigl\{0,\;\frac{\mathrm{s}}{T}\sum_{t\in\mathcal T} D_{b,t}\Bigr\},\qquad \mathrm{s}\in[0,1].
\end{equation}
The rebate $\rho_b\ge 0$ is a regulated constant that the building retains ($\rho_b=0$ recovers full extraction), and it is what keeps every building at or above its outside option.
\begin{proposition}[Hour-wise Individual Rationality (IR)]\label{prop:hourIR}
Under the fee \eqref{eq:reg-body} with share rule \eqref{eq:rho-proportional}, every building $b\in\mathcal B$ that plays its prescribed response satisfies in every hour $t\in\mathcal T$
\begin{equation} \label{eq: IR-constraint}
\Pi^{\mathrm B}_{b,t} - \Phi_{b,t} \;=\; \Pi^{\mathrm{Alt}}_{b,t}+\rho_b \;\ge\; \Pi^{\mathrm{Alt}}_{b,t},
\end{equation}
with strict inequality when $\rho_b>0$, and over the horizon, with $G_b=\sum_{t\in\mathcal T}D_{b,t}$ its participation gain,
\begin{equation}\label{eq:IR-horizon}
\sum_{t\in\mathcal T}\bigl(\Pi^{\mathrm B}_{b,t}-\Phi_{b,t}\bigr)=\sum_{t\in\mathcal T}\Pi^{\mathrm{Alt}}_{b,t}+\mathrm{s}\max\{0,G_b\}.
\end{equation}
\end{proposition}
Proposition~\ref{prop:hourIR} shows that no building is ever worse off than on its own, in any hour, a guarantee that neither the per-capita fixed part of \cite{vogelsang1989twopart} nor the HRGV fee \cite{hesamzadeh2018simple}, which ignores the users' benefit, can give. A building whose gain is positive keeps the share $\mathrm s$ of it as its dividend, and a building whose gain is negative is made whole by the supplier.

\begin{proposition}[Revenue Adequacy] \label{prop: prop-budeget}
Let $W=\sum_{t\in\mathcal T}W_t$ with $W_t=\Pi^{\mathrm{Sup}}_t+\sum_{b\in\mathcal B}\Pi^{\mathrm B}_{b,t}$ be the social welfare and $\Delta W=W-\sum_{b\in\mathcal B}\sum_{t\in\mathcal T}\Pi^{\mathrm{Alt}}_{b,t}$ the surplus the network adds over the audited standalone. Under the fee \eqref{eq:reg-body} with share rule \eqref{eq:rho-proportional} and the prescribed responses, the supplier's settled payoff is
\begin{equation}\label{eq:RA-identity}
\sum_{t\in\mathcal T}\Pi_t^{\mathrm{Sup}}+\sum_{b\in\mathcal B}\sum_{t\in\mathcal T}\Phi_{b,t}=\Delta W-\mathrm{s}\sum_{b\in\mathcal B}\max\{0,G_b\},
\end{equation}
and it is nonnegative, i.e., the mechanism is revenue adequate, if and only if
\begin{equation}\label{eq:RA-cond}
W\ge\sum_{b\in\mathcal B}\sum_{t\in\mathcal T}\Pi^{\mathrm{Alt}}_{b,t}+T\sum_{b\in\mathcal B}\rho_b .
\end{equation}
\end{proposition}
Based on Proposition~\ref{prop: prop-budeget}, the supplier thus recovers its costs whenever the network creates more value than the buildings would have on their own. 
The condition fails when no temperature policy lets the shared assets beat the standalone cost. The mechanism then runs a deficit, which signals that network operation is not welfare-improving. Vogelsang \cite{vogelsang1989twopart} notes that his tariff turns the regulated firm into a cooperative of its users, the fixed part being the membership contribution and the returned profit the dividend. The same holds here. The fee is the membership contribution, and the rebate is the dividend, paid per building in proportion to its own gain rather than per capita, and through the fixed part rather than through a price discount.
At the ceiling $\mathrm{s}_{\max}=\Delta W/\sum_b\max\{0,G_b\}$ of \eqref{eq:RA-identity} the supplier's payoff is zero and the whole surplus $\Delta W$ goes to the buildings, which is the cooperative outcome, reached here without cooperative premises and, unlike Ramsey pricing, without moving the price away from the marginal value of heat. For a share below the ceiling, the supplier keeps the part of $\Delta W$ that is not rebated.

\subsection{Resulting bilevel problem}\label{sec:bilevel}
This subsection combines the three-player problem. The regulatory fee $\Phi_{b,t}$ is determined solely by the regulatory constraints~\eqref{eq:reg-body} with no objective function, and is therefore merged into the supplier's constraints~\cite{hooman}. The result is the bilevel program below, whose lower level returns $(H^{\mathrm{Th,B}}_{b,t},\Pi^{\mathrm{B}}_{b,t})$ to the upper level:
\begin{flalign}\label{eq:UL-total}
&\begin{aligned}
\max_{\Omega^{\mathrm{UL}}}\ \ & \sum_{t\in\mathcal T}\!\Big(\Pi_t^{\mathrm{Sup}}+\!\sum_{b\in\mathcal B}\!\Phi_{b,t}\Big)-\upeta\sum_{b,t}\Phi^2_{b,t}\\
\mathrm{s.t.}\ \ & \text{regulatory constraints \eqref{eq:HRGV_phi_init}--\eqref{eq:HRGV_D},}\\
& \text{thermal supplier constraints \eqref{eq:UL-a}--\eqref{eq:UL-t},}\\
& \bigl(H^{\mathrm{Th,B}}_{b,t},\Pi^{\mathrm{B}}_{b,t}\bigr)\ \text{solve \eqref{eq:LL-total}},\quad\forall b\in\mathcal B,
\end{aligned}&&
\end{flalign}
\begin{flalign}\label{eq:LL-total}
&\hspace{2.6em}\begin{aligned}
\max_{\Omega^{\mathrm{LL}}}\ \ & \sum_{t\in\mathcal T}\bigl(\Pi^{\mathrm{B}}_{b,t}-\Phi_{b,t}\bigr)\\
\mathrm{s.t.}\ \ & \text{payoff and comfort definitions \eqref{eq:HRGV_Pi_B}--\eqref{eq:HRGV_x},}\\
& \text{building constraints \eqref{LL},}
\end{aligned}&&
\end{flalign}
where the lower level \eqref{eq:LL-total} is solved by each building for the posted upper-level decisions $(\lambda^{\mathrm{Th}}_t,\theta_{k,t})$, which enter it as parameters. Problem \eqref{eq:UL-total} is a mixed-integer nonlinear program,
with the binaries $\theta_{k,t},y^{\mathrm{inc}}_{t},y^{\mathrm{dec}}_{t}\in\{0,1\}$, the free tariff $\lambda^{\mathrm{Th}}_{t}\in\mathbb R$, the free net exchange, warming step and accumulator charge, all other variables continuous and nonnegative, and a bilinear revenue term.
Problem \eqref{eq:LL-total} is a convex quadratic program in continuous variables, nonnegative except $H^{\mathrm{Th,B}}_{b,t}$. The tariff revenue $\lambda^{\mathrm{Th}}_t\sum_b H^{\mathrm{Th,B}}_{b,t}$ appears with a positive sign in the supplier's profit \eqref{eq:UL-profit-ates-pen} and with a negative sign in the fee through \eqref{eq:HRGV_D}, so it cancels in the leader's objective. Welfare therefore pins the dispatch but not the price that supports it, i.e., the leader's problem is degenerate in $\lambda^{\mathrm{Th}}_t$, every welfare-supporting tariff giving the supplier the same payoff. This indifference leaves the buildings' settlement unaffected, since each ends at $\Pi^{\mathrm{Alt}}_{b,t}+\rho_b$ whatever supporting tariff is posted, but it would leave the price signal arbitrary. The degeneracy is resolved by adding a small fee-norm penalty $-\upeta\sum_{b,t}\Phi^2_{b,t}$ to the objective, which selects the tariff with minimum fee among the welfare optima.

To solve the bilevel program, the lower-level problem, convex for any fixed upper-level decision, is replaced by its Karush--Kuhn--Tucker (KKT) conditions in the upper level~\cite{hooman}, and the bilinear revenue term $\lambda^{\mathrm{Th}}_t\sum_b H^{\mathrm{Th,B}}_{b,t}$ in \eqref{eq:UL-profit-ates-pen} is replaced exactly by the lower-level strong-duality equality. The fee \eqref{eq:HRGV_D} inherits the squared comfort deviations of \eqref{eq:HRGV_utility}, so the penalty $-\upeta\sum_{b,t}\Phi^{2}_{b,t}$ would be a polynomial of degree four in the served levels. To keep the program convex, each squared deviation inside the penalized fee is replaced by its piecewise-linear interpolant \cite{Wu2011pwl}, which is written with a special ordered set of type 2 (SOS2) \cite{Keha2004}. So the single-level reformulation is a mixed-integer quadratic program. The complementarity conditions of the dual variables are linearized with big-$M$ constants. The KKT conditions and their related equations are provided in the online companion~\cite{Zenodo17362179}.

\subsection{Proposed tariff mechanism}\label{sec:protocol}
Now we are in a position to present the proposed mechanism, which is the announcement--dispatch--settlement protocol summarized in Algorithm~\ref{alg:regulatory_protocol}. First, the regulator announces the fee rule \eqref{eq:reg-body} and the sharing coefficient $\mathrm{s}$. Next, as the core, the utility solves the bilevel program \eqref{eq:UL-total}--\eqref{eq:LL-total}, following the approach explained in Section~\ref{sec:bilevel}, and then posts the tariff and the temperature schedule. Subsequently, the buildings self-dispatch against that posting and are billed at it from the meters, while the fee and the rebate, computed from \eqref{eq:reg-body}--\eqref{eq:rho-proportional} at the announced gains, stay fixed and are audited from the meters. In practice, Phase~4 settles a fee that was fixed before the day starts, so each building simply dispatches against the posted price, the response the supplier anticipated.
\begin{algorithm}[t]
\caption{Proposed Regulatory Mechanism: Announcement, Dispatch, and Settlement}
\label{alg:regulatory_protocol}
\begin{algorithmic}[1]

\Phase{Phase 1: Initialization}
\Statex \textbf{Data Ownership:}
\Statex $\bullet$ \textbf{Regulator:} Holds audited outside-alternatives, baseline contracts $(\hat{\mathrm{H}}^{h,c}_{b,t},\mathrm{d}^{\mathrm{h,c}}_b)$, and public tariffs.
\Statex $\bullet$ \textbf{Utility/Supplier:} Holds private asset data, cost functions, network models, and operational forecasts.
\Statex $\bullet$ \textbf{End Users/Buildings:} Hold private appliance constraints, domestic device profiles, and demand data.
\State Regulator announces the fee rule \eqref{eq:reg-body} and the gain-sharing coefficient $\mathrm{s}\in[0,1]$.

\Phase{Phase 2: Announcement}
\State Utility solves the bilevel program \eqref{eq:UL-total}--\eqref{eq:LL-total} and posts $\{\lambda^{\mathrm{Th}}_t, \theta_{k,t}\}$.
\State Regulator sets $\rho_b$ by \eqref{eq:rho-proportional} evaluated at the announced gains $\hat D_{b,t}$, $\forall b\in\mathcal B$; $\{\rho_b\}$ remains constant thereafter.

\Phase{Phase 3: Dispatch}
\For{$t=1,\dots,T$}
    \State Each building, in parallel, dispatches by \eqref{eq:LL-total} given $(\lambda^{\mathrm{Th}}_t,\theta_{k,t})$; sub-meters log $H^{\mathrm{Th,B}}_{b,t}$ and $P^{\mathrm{El,imp}}_{b,t}$.
\EndFor

\Phase{Phase 4: Settlement}
\For{each $b\in\mathcal B$ and $t\in\mathcal T$}
    \State $\Phi_{b,t} \leftarrow \hat D_{b,t} - \rho_b$, fixed at announcement; meters bill $\lambda^{\mathrm{Th}}_tH^{\mathrm{Th,B}}_{b,t}$.
\EndFor
\State Regulator audits $\hat D_{b,t}$ against the meters on demand.

\end{algorithmic}
\end{algorithm} 

\begin{proposition}[Static Efficiency] \label{prop: efficiency-prop}
Let $\overline W$ be the maximum of the social welfare $W$ of Proposition~\ref{prop: prop-budeget} over the supplier and building constraints of \eqref{eq:UL-total}--\eqref{eq:LL-total} jointly. If the supplier takes the rebates $\rho_b$ of \eqref{eq:rho-proportional} as given, then every global optimum of \eqref{eq:UL-total}--\eqref{eq:LL-total} with $\upeta=0$ attains $\overline W$, and for $\upeta>0$ the welfare loss satisfies $0\le\overline W-W\le\upeta\sum_{b\in\mathcal B}\sum_{t\in\mathcal T}(\Phi^{*}_{b,t})^{2}$, with $\Phi^{*}_{b,t}$ the fee at a welfare-optimal dispatch.\end{proposition}
As shown in Proposition~\ref{prop: efficiency-prop}, the fee hands the supplier the whole gain it creates for the buildings, so the schedule that maximizes its profit is the one that maximizes welfare. Under the lagged fixed part of \cite{vogelsang1989twopart} the firm earns only an approximation of each period's welfare change, whereas here the gain is exact and the optimum is reached at announcement. The share $\mathrm s$ therefore moves only the split of the gains, whereas the fee-tuning parameter of \cite{xia2024incentive} shares the surplus at the cost of efficiency. The fee-norm penalty is only the tie-breaker introduced above.

\begin{proposition}[Restricted-Marginal Pricing]\label{prop:price}
Let $\lambda^{*}_t$ be the shadow price of exchanged heat at the welfare-optimal dispatch, i.e., the multiplier of the loop balance \eqref{eq:UL-b} in the integer-restricted welfare program, shared by the exchange balances \eqref{eq:LL-res} of all buildings. The supplier takes the rebates $\rho_b$ of \eqref{eq:rho-proportional} as given. Then, at every hour in which some building responds to price at the margin, every profit-maximizing tariff satisfies $\lambda^{\mathrm{Th}}_t=\lambda^{*}_t$ and the supplier strictly prefers it.\end{proposition}
Based on Proposition \ref{prop:price}, the supplier posts the marginal value of heat, which the adjustment process of \cite{vogelsang1989twopart} reaches only in its steady state. Posting another price cannot increase the supplier's payoff, because the fee already transfers each building's whole gain to it. Also, a price away from $\lambda^{*}_t$ only lowers the welfare from which that gain is drawn, see \eqref{eq:app-harberger}, so no price cap is needed.

\section{Case Study}\label{sec:case_study}
This section evaluates the mechanism on an operating Belgian LTDHC, shown in Fig.~\ref{fig:network}. It consists of a central supplier and four office buildings on a bidirectional loop, with metered demand from a research campus in Genk (three office-hours loads, one continuously cooled facility) and Aquifer Thermal Energy Storage (ATES) parameterized from~\cite{rostam}, with a warm pipe (10--25\,\textdegree C) and a cold pipe (5--20\,\textdegree C) discretized into four temperature bins for piecewise COP interpolation~\cite{wirtz2021temperature}. On-network buildings use WSHPs (heating delivered at 60\,\textdegree C), WSChs (cooling at 16\,\textdegree C), and direct cooling from the cold pipe when the pinch condition is met. Rooftop PV output, precomputed from irradiance data at 18\% panel efficiency, enters settlement as the forecast $\hat{\mathrm P}^{\mathrm{PV}}_{b,t}$. An Air-Source Heat Pump (ASHP), a 92\%-efficient gas boiler, and an Air-Source Chiller (ASCh) appear only in the audited outside alternative $\Pi^{\mathrm{Alt}}_{b,t}$. The supplier operates a GSHP (warm well at 17\,\textdegree C, rated 180~kW) and GSCh (cold well at 6\,\textdegree C, rated 120~kW) with bin-wise COP curves, with the wells being at fixed temperature over the day,\footnote{Over one day the aquifer well temperatures are constant~\cite{rostam}, and the accumulator state and network temperature close cyclically so that the pricing scenarios are comparable.} and a 190~kWh\textsubscript{th} accumulator coupling the two pipes, so the architecture routes waste heat from cooling to heating users and the supplier covers only the net residual.
Demand is nearly inelastic, matching short-run district heating trends \cite{trotta2022elasticity}.
Hence, the comfort curvature is set to $\mathrm d=\bar{\mathrm p}/(2\upvarepsilon\bar{\mathrm H})$ \cite{Aalami2010}, with $\upvarepsilon=0.01$, $\bar{\mathrm p}$ the mean unit cost of thermal service, and $\bar{\mathrm H}$ the mean hourly load, so that a \(1\%\) deviation from the typical load costs the mean thermal service price. The utility \eqref{eq:HRGV_utility} carries no intercept, so the welfare levels reported below are negative and only differences across designs are meaningful.
Representative days are selected from the annual dataset by ranking every day by its heating share $\mathrm{h}_d=\sum_t\hat{\mathrm{H}}^{h}/\sum_t(\hat{\mathrm{H}}^{h}+\hat{\mathrm{H}}^{c})$, the fraction of the day's total thermal demand that is heating. The heating extreme is the day with the highest share ($\mathrm{h}_d\!=\!0.95$, winter Day of Year (DOY)~39), the cooling extreme is the day with the lowest share ($\mathrm{h}_d\!=\!0.04$, summer DOY~189), the balanced day is the one closest to an even split ($\mathrm{h}_d\!=\!0.49$, spring DOY~102), and the mixed autumn day matches Belgium's heating-dominated norm~\cite{eurostat_hdd} ($\mathrm{h}_d\!=\!0.63$, DOY~293). Electricity tariffs are time-of-use, 22.0 (off-peak) and 23.2 (on-peak)~c\euro/kWh for import and 4.0 and 7.0~c\euro/kWh for export, and gas, used only by the boiler of the outside option, costs 8.0~c\euro/kWh. The fee-norm weight is set to $\upeta=10^{-5}$, and every scenario is solved with Gurobi~13.0.2. Since ad hoc big-$M$ choices can truncate the dual space~\cite{PinedaMorales2019}, the big-$M$ constants are calibrated from the data and verified ex post. At the reported solutions, no dual variable approaches its bound, and the lower level, re-solved at the posted prices, reproduces the reported dispatch. Emissions are accounted for at a flat grid intensity $\mathrm{EF}=0.145$~kg/kWh on the supplier draw and the buildings' grid import, with PV exports not credited. The device ratings and network parameters are in the companion~\cite{Zenodo17362179}.

\subsection{Pricing scenarios}\label{ssec:pricing_scenarios}
The proposed mechanism is compared with three scenarios, S0 to S2, that share its physical model, data, and horizon but differ in the pricing rule and in who sets the price. All four scenarios are measured against the same standalone baseline $\Pi^{\mathrm{Alt}}_{b,t}$, the hourly standalone payoff per building, assuming no connection with the network and serving the same reference demand with an ASHP, a gas boiler, and an ASCh, with PV self-consumption or export, formulated in Appendix~\ref{app:alt}:
\begin{itemize}
\item \textbf{S0 (centralized welfare benchmark).}
A social planner maximizes the welfare $W$ of Proposition~\ref{prop: prop-budeget} over $\Omega^{\mathrm{CEN}}=(\Omega^{\mathrm{UL}}\cup\Omega^{\mathrm{LL}})\setminus\{\lambda^{\mathrm{Th}}_{t},\Phi_{b,t}\}$, subject to the thermal supplier constraints \eqref{eq:UL-a}--\eqref{eq:UL-t} and the building constraints \eqref{eq:HRGV_utility}--\eqref{eq:HRGV_x}, \eqref{LL}. The hourly value of heat is recovered ex-post as the dual of \eqref{eq:LL-res} in the fixed-integer convex relaxation. 
\item \textbf{S1 (Cap-based regulation).}
S1 is the same bilevel Stackelberg game structure as the proposed model, but the supplier maximizes its own profit $\sum_{t}\Pi^{\mathrm{Sup}}_t$ subject to the thermal supplier constraints \eqref{eq:UL-a}--\eqref{eq:UL-t}, the consumers' lower-level problem \eqref{eq:LL-total} with $\Phi_{b,t}\equiv0$, and the price limit constraints $-\mathrm{R}\,\overline{\uplambda}^C\le \lambda_t^{\mathrm{Th}}\le \mathrm{R}\,\overline{\uplambda}^H$ for every $t$, as in the capped variable heat tariff of \cite{Bai22}.
The cap references are the supplier's per-unit thermal cost at the peak electricity tariff and the best COP,
$\overline{\uplambda}^H=\uplambda^{\mathrm{El,imp}}_{\max}/\mathrm{COP}^{\mathrm{GSHP}}_{\max}=3.31$ and
$\overline{\uplambda}^C=\uplambda^{\mathrm{El,imp}}_{\max}/\mathrm{COP}^{\mathrm{GSCh}}_{\max}=3.86$~c\euro/kWh, and the headroom factor $\mathrm{R}=1.5$ leaves the supplier a $50\%$ margin above these costs in each direction of exchange~\cite{Lu21,Jin21,Mon24}. 
\item \textbf{S2 (Ramsey).}
S2 maximizes the welfare $W$ over $\Omega^{\mathrm{CEN}}\cup\{\lambda_t^{\mathrm{Th}}\}$, with the building responses embedded through the KKT conditions of \eqref{eq:LL-total} and an additional break-even constraint, $\sum_{b,t}\lambda_t^{\mathrm{Th}}H_{b,t}^{\mathrm{Th,B}} \ge \sum_{t}\uplambda_t^{\mathrm{El,imp}}P_t^{\mathrm{Sup}}$.
\item \textbf{Proposed mechanism.}
Algorithm~\ref{alg:regulatory_protocol} with the surplus share $\mathrm{s}=0.10$ in \eqref{eq:rho-proportional}, i.e., each building keeping a tenth of its gain over the audited standalone. 
\end{itemize}

\subsection{Results and Discussion}
\begin{figure}[t]
\centering
\includegraphics[width=0.9\linewidth]{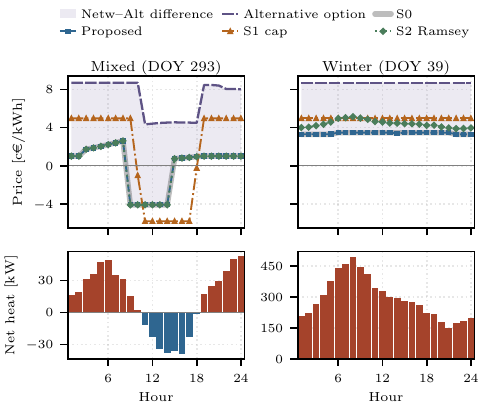}
\caption{Hourly tariffs (top) and net residual demand (bottom).}
\label{fig:tariffs}
\end{figure}

\begin{figure}[t]
\centering
\includegraphics[width=0.9\linewidth]{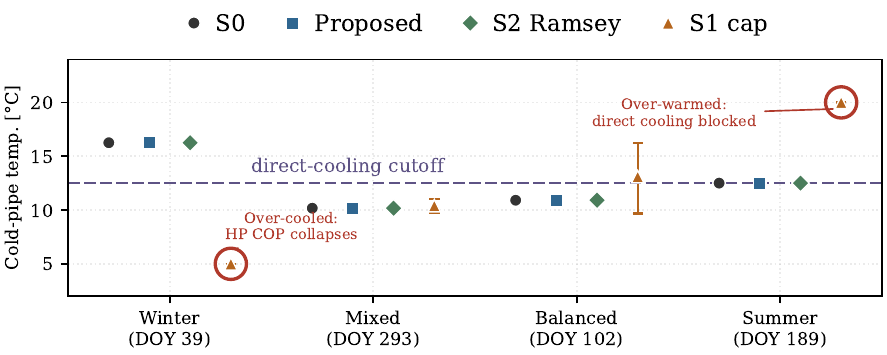}
\caption{Daily mean cold-pipe temperature per design (bars: S1 hourly range).}
\label{fig:network_dsp}
\end{figure}

\begin{figure}[t]
\centering
\includegraphics[width=0.9\linewidth]{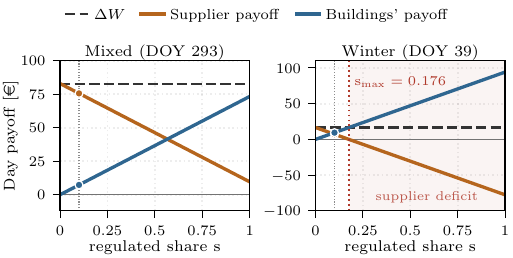}
\caption{Supplier and buildings' payoffs against the share $\mathrm{s}$ (dots: $\mathrm{s}=0.10$).}
\label{fig:sharing}
\end{figure}

Table~\ref{tab:welfare} gives the outcomes of the representative days. On all three days, the proposed mechanism reaches exactly the welfare of the centralized benchmark S0, as proven in Proposition~\ref{prop: efficiency-prop}. Under the proposed mechanism, on the mixed day, the supplier earns $75.69$~\euro\ and the buildings keep the $7.34$~\euro\ rebate. On the winter day, S1 turns the $16.63$~\euro\ surplus of the network into a $101.19$~\euro\ loss relative to standalone operation. S2 loses $0.13$~\euro, because its break-even binds and pushes the price $1.7$~c\euro/kWh above the shadow price, so the buildings consume slightly less than is efficient, to make the supplier break even. The proposed mechanism again matches S0 exactly, and its supplier still earns $7.21$~\euro\ (Proposition~\ref{prop: prop-budeget}), because the fee covers the loss on heat sales at the efficient price. On the summer day, the standalone option beats the network under every design, S0 included. S1 and S2 still operate it and pass the loss to the buildings, S1 even with a profit of $236.90$~\euro. Under the proposed mechanism the buildings are made whole. The loss falls on the supplier, so it declines to operate when operation destroys value.
\begin{table*}[t]\centering\footnotesize\setlength{\tabcolsep}{2.4pt}\renewcommand{\arraystretch}{0.60}
\caption{Day-level outcomes at $\mathrm{s}=0.10$. $\Delta W$ is the incremental surplus of Proposition~\ref{prop: prop-budeget} and $\Delta W$ loss the shortfall from S0.}\label{tab:welfare}
\begin{tabular}{l rrrr rrrr rrrr}
\toprule
 & \multicolumn{4}{c}{Mixed day (DOY 293)} & \multicolumn{4}{c}{Heating-extreme winter (DOY 39)} & \multicolumn{4}{c}{Cooling-extreme summer (DOY 189)} \\
\cmidrule(lr){2-5}\cmidrule(lr){6-9}\cmidrule(lr){10-13}
 & S0 & S1 cap & S2 Ramsey & Proposed & S0 & S1 cap & S2 Ramsey & Proposed & S0 & S1 cap & S2 Ramsey & Proposed \\
\midrule
Social welfare [\euro]    & $-47.01$ & $-47.93$ & $-47.01$ & $-47.01$ & $-834.87$ & $-952.69$ & $-835.00$ & $-834.87$ & $-175.18$ & $-355.24$ & $-175.18$ & $-175.18$ \\
$\Delta W$ [\euro]       & $83.03$  & $82.10$  & $83.03$  & $83.03$  & $16.63$   & $-101.19$ & $16.50$   & $16.63$   & $-180.59$ & $-360.65$ & $-180.59$ & $-180.59$ \\
$\Delta W$ loss [\euro]  & --       & $0.93$   & $0.00$   & $0.00$   & --        & $117.83$  & $0.13$    & $0.00$    & --        & $180.06$  & $0.00$    & $0.00$ \\
Supplier payoff [\euro]   & --       & $34.10$ & $9.67$  & $75.69$ & --        & $93.03$  & $0.00$    & $7.21$   & --        & $236.90$ & $0.00$    & $-180.59$ \\
Buildings' payoff [\euro] & --       & $48.00$  & $73.36$  & $7.34$   & --        & $-194.22$ & $16.50$   & $9.43$    & --        & $-597.56$ & $-180.59$ & $0.00$ \\
CO$_2$ [kg]               & $79.03$  & $78.83$  & $79.03$  & $79.03$  & $528.06$  & $569.75$  & $527.59$  & $528.06$  & $250.08$  & $298.53$  & $250.04$  & $250.10$ \\
\bottomrule
\end{tabular}
\end{table*}
Fig.~\ref{fig:tariffs} shows the price signal. The proposed tariff matches the efficient dispatch shadow price every hour, as shown in Proposition~\ref{prop:price}. On the mixed day it turns negative at midday, when chillers reject heat into the loop, meaning that the supplier pays for it. S2 posts the same price while its break-even is slack, and in winter it is $1.7$~c\euro/kWh more. S1 instead sits at its price cap and charges both directions of exchange.
Fig.~\ref{fig:network_dsp} shows the loop temperature, the control a price cap cannot reach. The proposed mechanism keeps it at the first-best on all four days (as do S0 and S2), because the fee makes the supplier internalize the buildings' efficiency. S1 moves it for profit: $5$\,\textdegree C in winter, lowering the heat-pump COP, and $20$\,\textdegree C on the cold pipe in summer, disabling direct cooling. The buildings lose more than the supplier gains, the difference dissipated in their equipment, as price regulation gives weak quality incentives \cite{WEISMAN2005165}. Seasonal storage in the warm well~\cite{rostam} is not counted, so the summer verdict covers the day only.
Table~\ref{tab:perbldg} breaks the settlement down by building. Every gain $\sum_tD_{b,t}$ is positive, and the fee is that gain less the rebate $T\rho_b$, so each building keeps a tenth of its gain above its outside option, as proven in Proposition~\ref{prop:hourIR}. Under S1, three buildings end $41$--$130$~\euro\ below their standalone payoff in winter, which Proposition~\ref{prop:hourIR} rules out for the proposed mechanism. Office~2, the cooling-only building, is paid for the heat it rejects under both designs, but gains $33$~\euro\ under S1 instead of $18$~\euro, because the $5$\,\textdegree C loop lets it cool directly without running its chiller and the cap price pays more for its heat.
\begin{table}[t]\centering\footnotesize\setlength{\tabcolsep}{3pt}\renewcommand{\arraystretch}{0.55}
\caption{Per-building outcomes [\euro] at $\mathrm{s}=0.10$.}\label{tab:perbldg}
\begin{tabular}{lrrrrr}
\toprule
Building & $\Pi^{\mathrm{Alt}}_b$ & $\sum_tD_{b,t}$ & $\Phi_b$ & $T\rho_b$ & S1 gain \\
\midrule
\multicolumn{6}{l}{\textbf{Mixed day (DOY 293)}}\\
Office 1       & $-49.00$ & $18.08$ & $16.27$ & $1.81$  & $0.59$ \\
Office 2       & $6.51$   & $7.66$  & $6.89$  & $0.77$  & $9.49$ \\
Office 3       & $-38.77$ & $27.96$ & $25.16$ & $2.80$  & $26.44$ \\
Office 4       & $-48.77$ & $19.68$ & $17.71$ & $1.97$  & $11.49$ \\
\midrule
\multicolumn{6}{l}{\textbf{Heating-extreme winter day (DOY 39)}}\\
Office 1       & $-514.98$ & $45.48$ & $40.93$ & $4.55$  & $-130.06$ \\
Office 2       & $-9.27$   & $18.15$ & $16.34$ & $1.82$  & $33.01$ \\
Office 3       & $-163.47$ & $16.42$ & $14.77$ & $1.64$  & $-56.64$ \\
Office 4       & $-163.77$ & $14.21$ & $12.79$ & $1.42$  & $-40.53$ \\
\bottomrule
\end{tabular}
\end{table}
Finally, Fig.~\ref{fig:sharing} shows that the share $\mathrm{s}$ is purely distributional. It changes neither dispatch, tariff, nor welfare (Proposition~\ref{prop: efficiency-prop}), only the split of $\Delta W$, linearly by \eqref{eq:RA-identity}. The supplier stays solvent for any $\mathrm{s}\le1$ on the mixed day ($\mathrm{s}_{\max}=1.13$) but only up to $0.176$ in winter.

\section{Conclusions} \label{sec:conclusion}
This paper proposed an incentive-based regulatory mechanism for tariff design in LTDHC networks, anchoring the regulated fee to an audited outside alternative instead of the firm's history. Three findings stand out from the Belgian case study. First, the mechanism aligns both instruments of the supplier, the tariff and the loop temperature, and declines to operate the LTDHC when operation destroys value. Second, the proposed mechanism recovers costs without changing the efficient price and without charging any user more than its outside option. Third, the regulated share is a purely distributional dial without any welfare loss. Future work extends the mechanism to a stochastic version that includes forecast uncertainty, a decomposition algorithm, and to building-specific sharing rules.

\section*{Acknowledgment}
This research was carried out in the context of the Energy Transition Fund, REINVENT project, supported by the FPS Economy, S.M.E.s, Self-Employed and Energy, Belgium. The authors thank VITO--EnergyVille for the campus
metering data and the network specifications underlying the case study.

\appendices
\section{Proofs}\label{app:proofs}

\begin{IEEEproof}[Proof of Proposition~\ref{prop:brc}]
Fix the announced $(\lambda^{\mathrm{Th}}_t,\theta_{k,t})$. The fee is fixed at announcement and enters \eqref{eq:LL-total} as a constant, $\partial\Phi^{*}_{b,t}/\partial\Omega^{\mathrm{LL}}_b=0$, so $\arg\max_{\Omega^{\mathrm{LL}}_b}(\Pi^{\mathrm{B}}_{b,t}-\Phi^{*}_{b,t})=\arg\max_{\Omega^{\mathrm{LL}}_b}\Pi^{\mathrm{B}}_{b,t}$. The lower level maximizes a concave quadratic over linear constraints and satisfies Slater's condition, so its KKT conditions are necessary and sufficient, and the strictly concave comfort function \eqref{eq:HRGV_utility} makes the served heat of the maximizer unique. Any unilateral deviation that changes the served heat or raises the bill lowers $\Pi^{\mathrm{B}}_{b,t}$ and leaves $\Phi^{*}_{b,t}$ unchanged, hence lowers the net payoff.
\end{IEEEproof}

\begin{IEEEproof}[Proof of Proposition~\ref{prop:info}]
Both payoffs in \eqref{eq:HRGV_D} are the contracted comfort function \eqref{eq:HRGV_utility} minus a purchase cost, so $D_{b,t}=[U_{b,t}(x^{\mathrm B})-U_{b,t}(x^{\mathrm{Alt}})]-[\mathrm C^{\mathrm B}_{b,t}-\mathrm C^{\mathrm{Alt}}_{b,t}]$. With inelastic demand $x^{\mathrm B}=x^{\mathrm{Alt}}=\hat{\mathrm H}$, the first bracket vanishes and $D_{b,t}=\mathrm C^{\mathrm{Alt}}_{b,t}-\mathrm C^{\mathrm B}_{b,t}$. With elastic demand let $\mathrm{m}^{\mathrm B}$ and $\mathrm{m}^{\mathrm{Alt}}$ be the marginal costs of one served unit on and off the network, $\bar{\mathrm m}$ their mean and $\Delta\mathrm m=\mathrm{m}^{\mathrm{Alt}}-\mathrm{m}^{\mathrm B}$. Each side maximizes $-\mathrm d(x-\hat{\mathrm H})^{2}$ minus its purchase cost, so at an interior optimum $x-\hat{\mathrm H}=-\mathrm m/2\mathrm d$ and the comfort term equals $-\mathrm m^{2}/4\mathrm d$. Hence $U(x^{\mathrm B})-U(x^{\mathrm{Alt}})=[(\mathrm{m}^{\mathrm{Alt}})^{2}-(\mathrm{m}^{\mathrm B})^{2}]/4\mathrm d=\Delta\mathrm m\,\bar{\mathrm m}/2\mathrm d$ and $D_{b,t}=\mathrm C^{\mathrm{Alt}}_{b,t}-\mathrm C^{\mathrm B}_{b,t}+\Delta\mathrm m\,\bar{\mathrm m}/2\mathrm d$. With the curvature calibrated to the elasticity, $\mathrm d=\bar{\mathrm p}/(2\upvarepsilon\bar{\mathrm H})$ for a mean unit cost $\bar{\mathrm p}$ and a mean load $\bar{\mathrm H}$, the residual equals $\upvarepsilon(\bar{\mathrm m}/\bar{\mathrm p})\bar{\mathrm H}\Delta\mathrm m=\mathcal O(\upvarepsilon)$, computed from $\mathrm d$ and the served levels \eqref{eq:HRGV_x}.
\end{IEEEproof}

\begin{IEEEproof}[Proof of Proposition~\ref{prop:hourIR}]
At the prescribed response, by \eqref{eq:HRGV_phi_init}--\eqref{eq:HRGV_D}, $\Phi_{b,t}=\Phi_{b,t-1}+D_{b,t}-D_{b,t-1}=D_{b,t}-\rho_b$ and $D_{b,t}=\Pi^{\mathrm{B}}_{b,t}-\Pi^{\mathrm{Alt}}_{b,t}$, so for every $t\in\mathcal T$,
\begin{equation}\label{eq:phi-telescoped}
\Phi_{b,t}=D_{b,t}-\rho_b=\Pi^{\mathrm{B}}_{b,t}-\bigl(\Pi^{\mathrm{Alt}}_{b,t}+\rho_b\bigr).
\end{equation}
Rearranging gives $\Pi^{\mathrm B}_{b,t}-\Phi_{b,t}=\Pi^{\mathrm{Alt}}_{b,t}+\rho_b$ for every $b$ and $t$, the equality in \eqref{eq: IR-constraint}. The rebate is nonnegative by \eqref{eq:rho-proportional}, which gives the inequality, strict when $\rho_b>0$. Summing the equality over $t\in\mathcal T$ gives $\sum_t(\Pi^{\mathrm B}_{b,t}-\Phi_{b,t})=\sum_t\Pi^{\mathrm{Alt}}_{b,t}+T\rho_b$, and inserting $T\rho_b=\mathrm{s}\max\{0,G_b\}$ from \eqref{eq:rho-proportional} gives \eqref{eq:IR-horizon}. The fees of the day thus sum to $(1-\mathrm s)G_b$ for a positive gain and to $G_b$ itself for a negative one.
\end{IEEEproof}

\begin{IEEEproof}[Proof of Proposition~\ref{prop: prop-budeget}]
Substituting \eqref{eq:phi-telescoped} into the supplier's settled payoff \eqref{eq:UL-obj-short},
\begin{align}\label{eq: phi-with-alt}
\sum_{t\in\mathcal T}\Bigl(\Pi_t^{\mathrm{Sup}}+\sum_{b\in\mathcal B}\Phi_{b,t}\Bigr)
&=\sum_{t\in\mathcal T}\Bigl(\Pi_t^{\mathrm{Sup}}+\sum_{b\in\mathcal B}\Pi_{b,t}^{\mathrm B}\Bigr)-\mathrm{K}\nonumber\\
&=W-\mathrm K,
\end{align}
with $\mathrm K=\sum_{b,t}(\Pi^{\mathrm{Alt}}_{b,t}+\rho_b)=\sum_{b,t}\Pi^{\mathrm{Alt}}_{b,t}+T\sum_b\rho_b$ a constant, with audited payoffs and fixed rebates. Inserting $T\rho_b=\mathrm{s}\max\{0,G_b\}$ from \eqref{eq:rho-proportional} gives $W-\sum_{b,t}\Pi^{\mathrm{Alt}}_{b,t}-\mathrm s\sum_b\max\{0,G_b\}=\Delta W-\mathrm{s}\sum_{b}\max\{0,G_b\}$, which is \eqref{eq:RA-identity}. Since $T\sum_b\rho_b=\mathrm s\sum_b\max\{0,G_b\}$, the right-hand side of \eqref{eq:RA-identity} is nonnegative exactly when \eqref{eq:RA-cond} holds.
\end{IEEEproof}

\begin{IEEEproof}[Proof of Proposition~\ref{prop: efficiency-prop}]
Set $\upeta=0$. By \eqref{eq: phi-with-alt} the supplier maximizes $W-\mathrm K$ with $\mathrm K$ a constant, so a global optimum of \eqref{eq:UL-total}--\eqref{eq:LL-total} attains $\overline W$, which is feasible for \eqref{eq:UL-total} by Proposition~\ref{prop:price}.

For $\upeta>0$ let $z$ denote a dispatch of the plant and the buildings and $F(z)=W(z)-\mathrm K-\upeta\sum_{b,t}\Phi_{b,t}(z)^{2}$ the leader's objective, with fees from \eqref{eq:phi-telescoped} at $z$. Let $z_{\upeta}$ maximize $F$ and $z^{*}$ maximize $W$. Since $z^{*}$ maximizes $W$, $W(z^{*})\ge W(z_{\upeta})$, and since $z_{\upeta}$ maximizes $F$, $F(z_{\upeta})\ge F(z^{*})$, i.e.,
\begin{align}
&W(z^{*})-\upeta\!\sum_{b,t}\!\bigl(\Phi^{*}_{b,t}\bigr)^{2}\ \le\ W(z_{\upeta})-\upeta\!\sum_{b,t}\!\Phi_{b,t}(z_{\upeta})^{2}.\label{eq:eta-opt}
\end{align}
Rearranging and dropping $\upeta\sum_{b,t}\Phi_{b,t}(z_{\upeta})^{2}\ge0$,
\begin{align}\label{eq:eta-bound}
0&\le W(z^{*})-W(z_{\upeta})\nonumber\\
&\le\upeta\sum_{b,t}\Bigl[\bigl(\Phi^{*}_{b,t}\bigr)^{2}-\Phi_{b,t}(z_{\upeta})^{2}\Bigr]
\le\upeta\sum_{b,t}\bigl(\Phi^{*}_{b,t}\bigr)^{2}.
\end{align}
\end{IEEEproof}

\begin{IEEEproof}[Proof of Proposition~\ref{prop:price}]
Solve the welfare program, the maximization of $W$ over the constraints of \eqref{eq:UL-total}--\eqref{eq:LL-total}, once as a mixed-integer program and fix its integer variables $(\theta_{k,t},y^{\mathrm{inc}}_{t},y^{\mathrm{dec}}_{t})$ at their optimal values. The remaining program is convex with the same optimal value, and its multipliers are supporting prices, the discrete residual being settled by contracts, here the fee \eqref{eq:reg-body} \cite[Theorems~1--2]{ONEILL2005269}. Let $z^{*}$ be its optimal dispatch of the plant and of all buildings, $\overline W$ its welfare, $\lambda^{*}_t$ the multiplier of the loop balance \eqref{eq:UL-b} at $z^{*}$, and $z(\lambda^{\mathrm{Th}})$ the buildings' best response to a posted tariff (Proposition~\ref{prop:brc}). The exchange $H^{\mathrm{Th,B}}_{b,t}$ is free and appears in \eqref{eq:UL-b} and in \eqref{eq:LL-res}, so stationarity in it sets the multiplier of every \eqref{eq:LL-res} equal to $\lambda^{*}_t$.

Let $\mathcal L$ be the Lagrangian of \eqref{eq:LL-obj}--\eqref{LL}. Stationarity in the active-bin heat-pump power, $\partial\mathcal L/\partial P^{\theta,\mathrm{Hp}}_{b,k,t}=0$, gives
\begin{equation}\label{eq:app-stat}
\mathrm{COP}_k\,2\mathrm d\,(x-\hat{\mathrm H})+(\mathrm{COP}_k-1)\,\nu_t+\delta_t+\mu^{h}-\nu^{\theta\mathrm{Hp}}=0,
\end{equation}
with $x$ the served heat \eqref{eq:HRGV_x}, $\hat{\mathrm H}\equiv\hat{\mathrm H}^{h}_{b,t}$, $\mathrm d\equiv\mathrm d^{\mathrm h}_b$, $\nu_t$ and $\delta_t$ the duals of \eqref{eq:LL-res} and \eqref{eq:LL-Pbal}, $\mu^{h}$ the dual of the bin cap in \eqref{eq:LL-cap} and $\nu^{\theta\mathrm{Hp}}$ the dual of its nonnegativity bound. For an interior response $\mu^{h}=\nu^{\theta\mathrm{Hp}}=0$, stationarity in $H^{\mathrm{Th,B}}_{b,t}$ gives $\nu_t=\lambda^{\mathrm{Th}}_t$, and stationarity in $P^{\mathrm{El,imp}}_{b,t}$ gives $\delta_t=\uplambda^{\mathrm{El,imp}}_t$ at import hours (export hours change only $\kappa$ below). Dividing \eqref{eq:app-stat} by $\mathrm{COP}_k$, with $c:=(\mathrm{COP}_k-1)/\mathrm{COP}_k$ the metering constant $\partial H^{\mathrm{Th,B}}/\partial x$ and $\kappa:=\uplambda^{\mathrm{El,imp}}_t/\mathrm{COP}_k$ the electricity cost per served unit,
\begin{equation}\label{eq:app-foc}
2\mathrm d\bigl(\hat{\mathrm H}-x\bigr)=c\,\lambda^{\mathrm{Th}}_t+\kappa
\;\Longrightarrow\;
x(\lambda^{\mathrm{Th}}_t)=\hat{\mathrm H}-\frac{c\,\lambda^{\mathrm{Th}}_t+\kappa}{2\mathrm d},
\end{equation}
with slope $dx/d\lambda^{\mathrm{Th}}_t=-c/2\mathrm d$ (cooling differs by $c=(\mathrm{COP}_k+1)/\mathrm{COP}_k$ and the price sign). Next, write the hourly welfare as a function of $x$. The building contributes the comfort term $-\mathrm d(x-\hat{\mathrm H})^{2}$ and the electricity cost $\kappa\,x$, the heat it draws from the loop, $c\,x$, is valued at $\lambda^{*}_t$, and every other term of $W_t$ is independent of $x$ and collected in $\mathrm K'$, so
\begin{equation}\label{eq:app-sw}
W_t(x)=-\mathrm d\bigl(x-\hat{\mathrm H}\bigr)^{2}-\kappa\,x-\lambda^{*}_t c\,x+\mathrm K'.
\end{equation}
Let $\mathrm m:=c\lambda^{\mathrm{Th}}_t+\kappa$ be the cost of one served unit to the building at the posted tariff and $\mathrm m^{*}:=c\lambda^{*}_t+\kappa$ its cost to the system, at the shadow price. By \eqref{eq:app-foc}, $x-\hat{\mathrm H}=-\mathrm m/2\mathrm d$, so the comfort term of \eqref{eq:app-sw} equals $-\mathrm m^{2}/4\mathrm d$ and the cost term $-\mathrm m^{*}x$ equals $\mathrm m^{*}\mathrm m/2\mathrm d$ up to a constant. Completing the square,
\begin{equation}\label{eq:app-square}
-\frac{\mathrm m^{2}}{4\mathrm d}+\frac{\mathrm m^{*}\mathrm m}{2\mathrm d}=-\frac{(\mathrm m-\mathrm m^{*})^{2}}{4\mathrm d}+\mathrm{const},
\end{equation}
with $\mathrm m-\mathrm m^{*}=c(\lambda^{\mathrm{Th}}_t-\lambda^{*}_t)$, and $\kappa$ cancels. Summing over the buildings whose response is interior,
\begin{equation}\label{eq:app-harberger}
W\bigl(z(\lambda^{\mathrm{Th}}_t)\bigr)\le\overline W-\sum_{b\,\text{active}}\frac{c_{b,t}^{2}}{4\mathrm{d}^{\mathrm{h}}_b}\bigl(\lambda^{\mathrm{Th}}_t-\lambda^{*}_t\bigr)^{2}.
\end{equation}
By \eqref{eq: phi-with-alt} the supplier's profit is $W-\mathrm K$. At an hour with a responding building, every tariff other than $\lambda^{*}_t$ therefore lowers welfare and profit.
\end{IEEEproof}

\section{Outside-Alternative Benchmark}\label{app:alt}
The audited outside option is the building's own problem off the network, which chooses the served levels $x^{\mathrm h},x^{\mathrm c}$ under the comfort preferences \eqref{eq:HRGV_utility} with an air-source heat pump, a gas boiler, and an air-source chiller.
For each building $b$ and hour $t$, over nonnegative served levels and powers:%
\begin{align*}
\Pi^{\mathrm{Alt}}_{b,t}=\max\ & U_{b,t}-\uplambda^{\mathrm{gas}}P^{\mathrm{GB}}_{b,t}-\uplambda^{\mathrm{El,imp}}_t P^{\mathrm{a,El}}_{b,t}+\uplambda^{\mathrm{El,exp}}_t P^{\mathrm{a,exp}}_{b,t}\\
\mathrm{s.t.}\ \ & \mathrm{COP}^{\mathrm{ASHP}}_t P^{\mathrm{ASHP}}_{b,t}+\upeta^{\mathrm{GB}}P^{\mathrm{GB}}_{b,t}=x^{\mathrm h},\\
& \mathrm{COP}^{\mathrm{ASCh}}_t P^{\mathrm{ASCh}}_{b,t}=x^{\mathrm c},\ \ \hat{\mathrm P}^{\mathrm{PV}}_{b,t}=P^{\mathrm{a,exp}}_{b,t}+P^{\mathrm{a,s}}_{b,t},\\
& P^{\mathrm{a,El}}_{b,t}=P^{\mathrm{ASHP}}_{b,t}+P^{\mathrm{ASCh}}_{b,t}-P^{\mathrm{a,s}}_{b,t},\\
& P^{\mathrm{ASHP}}_{b,t}\le\overline{\mathrm P}^{\mathrm{ASHP}}_b,\qquad P^{\mathrm{GB}}_{b,t}\le\overline{\mathrm P}^{\mathrm{GB}}_b,\\
& P^{\mathrm{ASCh}}_{b,t}\le\overline{\mathrm P}^{\mathrm{ASCh}}_b,
\end{align*}

\bibliographystyle{IEEEtran}
\bibliography{biblio}
\end{document}